\documentclass[12pt]{amsart}
\usepackage[mathscr]{eucal}
\usepackage{amsthm}
\usepackage{amssymb}
\usepackage{amscd}
\usepackage{verbatim}
\usepackage{pb-diagram}

\def\Bbb{\mathbb}
\def\frak{\mathfrak}

\newenvironment{pf*}[1]{\proof[#1]}{\endproof}
\newcommand{\rom}{\textup}
\newtheorem{Theorem}[equation]{Theorem}
\newtheorem{Corollary}[equation]{Corollary}
\newtheorem{Lemma}[equation]{Lemma}
\newtheorem{Proposition}[equation]{Proposition}
\newtheorem{Axiom}{Axiom}
\newtheorem{Conjecture}[equation]{Conjecture}

\theoremstyle{definition}
\newtheorem{Definition}[equation]{Definition}

\theoremstyle{remark}
\newtheorem{Remark}[equation]{Remark}

\numberwithin{equation}{section}

\newcommand{\thmref}[1]{Theorem~\ref{#1}}
\newcommand{\secref}[1]{\S\ref{#1}}
\newcommand{\lemref}[1]{Lemma~\ref{#1}}
\newcommand{\propref}[1]{Proposition~\ref{#1}}
\newcommand{\corref}[1]{Corollary~\ref{#1}}

\newcommand{\defeq}{\overset{\operatorname{\scriptstyle def.}}{=}}
\newcommand{\C}{{\Bbb C}}
\newcommand{\Z}{{\Bbb Z}}
\newcommand{\Q}{{\Bbb Q}}

\newcommand{\GL}{\operatorname{GL}}

\newcommand{\algsl}{\operatorname{\frak{sl}}} 

\newcommand{\g}{{\frak g}}

\newcommand{\Hom}{\operatorname{Hom}}

\newcommand{\Ker}{\operatorname{Ker}}

\newcommand{\Ima}{\operatorname{Im}}

\newcommand{\rank}{\operatorname{rank}}

\newcommand{\id}{\operatorname{id}}
\newcommand{\ve}{\varepsilon}
\newcommand{\vin}{\operatorname{in}} 
\newcommand{\vout}{\operatorname{out}} 
\newcommand{\M}{{\frak M}} 

\newcommand{\dslash}{/\!\!/} 
\newcommand{\bw}{{\mathbf w}} 

\newcommand{\codim}{\operatorname{codim}} 

\newcommand{\Ua}{{\mathbf U}_q(\widehat{\mathfrak g})} 
\newcommand{\Ul}{{\mathbf U}_q({\mathbf L}{\mathfrak g})} 

\newcommand{\Ue}{{\mathbf U}_{\varepsilon}(\mathfrak g)}

\newcommand{\Uli}{{\mathbf U}^{\Z}_q({\mathbf L}\mathfrak g)}
\newcommand{\Ule}{{\mathbf U}_{\varepsilon}({\mathbf L}\mathfrak g)}

\newcommand{\Zl}{\widetilde{\mathfrak Z}^\bullet}
\newcommand{\Zm}{\mathfrak Z^\bullet}
\newcommand{\Lg}{\mathbf L\g}
\newcommand{\bfR}{\mathbf R}
\newcommand{\N}{\mathfrak M^\bullet}
\newcommand{\Nreg}{\mathfrak M^{\bullet\operatorname{reg}}_0}
\newcommand{\NLa}{\mathfrak L^\bullet}
\newcommand{\bM}{{\mathbf M}^\bullet} 
\newcommand{\HomE}{\operatorname{E}^\bullet}
\newcommand{\HomL}{\operatorname{L}^\bullet}

\begin{document}
\title[$t$--analogs of $q$--characters]
{Quiver varieties and
$t$--analogs of $q$--characters of quantum affine algebras
}
\author{Hiraku Nakajima}
\address{Department of Mathematics, Kyoto University, Kyoto 606-8502,
Japan
}
\email{nakajima@kusm.kyoto-u.ac.jp}
\urladdr{http://www.kusm.kyoto-u.ac.jp/\textasciitilde nakajima}
\thanks{Supported by the Grant-in-aid
for Scientific Research (No.11740011), the Ministry of Education,
Japan.}
%
\begin{abstract}
Let us consider a specialization of an untwisted quantum affine
algebra of type $ADE$ at a nonzero complex number, which may or may
not be a root of unity.
The Grothendieck ring of its finite dimensional representations has
two bases, simple modules and standard modules.
We identify entries of the transition matrix with special values of
``computable'' polynomials, similar to Kazhdan-Lusztig polynomials.
At the same time we ``compute'' $q$-characters for all simple modules.
The result is based on ``computations'' of Betti numbers of
graded/cyclic quiver varieties.
(The reason why we put `` " will be explained in the end of the
introduction.)
\end{abstract}
\maketitle
\tableofcontents
\section*{Introduction}

Let $\g$ be a simple Lie algebra of type $ADE$ over $\C$, $\Lg = \g
\otimes \C[z,z^{-1}]$ be its loop algebra, and $\Ul$ be its quantum
universal enveloping algebra, or the quantum loop algebra for short.
It is a subquotient of the quantum affine algebra $\Ua$, i.e., without
central extension and degree operator. Let $\Ule$ be its
specialization at $q=\ve$, a nonzero complex number. (See
\secref{sec:qloop} for definition.)

It is known that $\Ule$ is a Hopf algebra. Therefore the category
$\mathscr Rep\Ule$ of finite dimensional representations of $\Ule$ is
a monoidal (or tensor) abelian category. Let $\operatorname{Rep}\Ule$
be its Grothendieck ring. It is known that $\operatorname{Rep}\Ule$ is
commutative (see e.g., \cite[Corollary~2]{FR}).

The ring $\operatorname{Rep}\Ule$ has two natural bases, simple
modules $L(P)$ and standard modules $M(P)$, where $P$ is the Drinfeld
polynomial. The latters were introduced by the author \cite{Na-qaff}.
The purpose of this article is to ``compute'' the transition matrix
between two bases.
(The reason why we put `` " will be explained in the end of the
introduction.)
More precisely, we define certain ``computable'' polynomials
$Z_{PQ}(t)$, which are analogs of Kazhdan-Lusztig polynomials for
Weyl groups. Then we show that the multiplicity $[M(P):L(Q)]$ is equal
to $Z_{PQ}(1)$.
This generalizes a result of Arakawa~\cite{Arakawa} who expressed the
multiplicities by Kazhdan-Lusztig polynomials when $\g$ is of type
$A_n$ and $\ve$ is not a root of unity.
Furthermore, coefficients of $Z_{PQ}(t)$ are equal to multiplicities of
simple modules of subquotients of standard module with respect to a
Jantzen filtration if we combine our result with \cite{Gr}, where the
transvesal slice used there was given in \cite{Na-qaff}.

Since there is a slight complication when $\ve$ is a root of unity, we 
assume $\ve$ is {\it not\/} so in this introduction.
Then the definition of $Z_{PQ}(t)$ is as follows.
Let $\bfR_t \defeq \operatorname{Rep}\Ule\otimes_\Z\Z[t, t^{-1}]$,
which is a $t$--analog of the representation ring.
By \cite{Na-qaff}, $\bfR_t$ is identified with the dual of the
Grothendieck group of a category of perverse sheaves on affine graded
quiver varieties (see \secref{sec:quiver} for the definition) so that
(1) $\{ M(P)\}$ is the specialization at $t=1$ of the dual base of
constant sheaves of strata, extended by $0$ to the complement,
and
(2) $\{ L(P)\}$ is that of the dual base of intersection cohomology
sheaves of strata.
%
%
A property of intersection cohomology complexes leads to the
following combinatorial definition of $Z_{PQ}(t)$:
Let $\setbox5=\hbox{A}\overline{\rule{0mm}{\ht5}\hspace*{\wd5}}$
be the involution on $\bfR_t$, dual to the Grothendieck-Verdier duality.
We denote the two bases of $\bfR_t$ by the same symbols $M(P)$, $L(P)$
at the specialization at $t=1$ for simplicity.
Let us express the involution in the basis $\{ M(P) \}_P$, classes of
standard modules:
\begin{equation*}
    \overline{M(P)} = \sum_{Q: Q\le P} u_{PQ}(t) M(Q),
\end{equation*}
where $\le$ is a certain ordering $<$ among $P$'s. We then define an
element $L(P)$ by
\begin{equation}\label{eq:can}
    \overline{L(P)} = L(P), \qquad
    L(P) \in M(P) + \sum_{Q: Q<P} t^{-1}\Z[t^{-1}] M(Q).
\end{equation}
The above polynomials $Z_{PQ}(t)\in \Z[t^{-1}]$ are given by
\begin{equation*}
    M(P) = \sum_{Q: Q\le P} Z_{PQ}(t) L(Q).
\end{equation*}
The existence and uniqueness of $L(P)$ (and hence of $Z_{PQ}(t)$) is
proved exactly as in the case of the Kazhdan-Lusztig polynomial.
In particular, it gives us a combinatorial algorithm computing
$Z_{PQ}(t)$, once $u_{PQ}(t)$ is given.

In summary, we have the following analogy:
\begin{center}
\begin{tabular}{|c|c|}
\hline
$\bfR_t$ & the Iwahori-Hecke algebra $H_q$
\\
standard modules $\{ M(P)\}_P$ & $\{ T_w \}_{w\in W}$
\\
simple modules $\{ L(P)\}_P$ & Kazhdan-Lusztig basis $\{ C'_w \}_{w\in W}$
\\
\hline
\end{tabular}
\end{center}
See \cite{KL} for definitions of $H_q$, $T_w$, 
$C'_w$.
%

The remaining task is to ``compute'' $u_{PQ}(t)$.
For this purpose we introduce a $t$--analog $\widehat\chi_{\ve,t}$
of the $q$--character, or $\ve$--character.
The original $\ve$--character $\chi_\ve$, which is a specialization of
our $t$--analog at $t=1$, was introduced by Knight~\cite{Kn} (for
Yangian and generic $\ve$) and Frenkel-Reshetikhin~\cite{FR} (for
generic $\ve$) and Frenkel-Mukhin~\cite{FM2} (when $\ve$ is a root of
unity).
It is an injective ring homomorphism from $\operatorname{Rep}\Ule$ to
$\Z[Y_{i,a}^\pm]_{i\in I,a\in\C^*}$, a ring of Laurent polynomials of
infinitely many variables.
It is an analog of the ordinary character homomorphism of the finite
dimensional Lie algebra $\mathfrak g$.
Our $t$--analog is an injective $\Z[t,t^{-1}]$-linear map
\begin{equation*}
   \widehat\chi_{\ve,t} \colon 
   \bfR_t
   \to \widehat{\mathscr Y}_t \defeq
       \Z[t, t^{-1}, V_{i,a}, W_{i,a}]_{i\in I,a\in\C^*}.
\end{equation*}
We have a simple, explicit definition of an involution
$\setbox5=\hbox{A}\overline{\rule{0mm}{\ht5}\hspace*{\wd5}}$ on
$\widehat{\mathscr Y}_t$ (see \eqref{eq:bar}). The involution on
$\bfR_t$ is the restriction. Therefore the matrix $(u_{PQ}(t))$ can be
expressed in terms of values of $\widehat\chi_{\ve,t}(M(P))$ for all
$P$.

We define $\widehat\chi_{\ve,t}$ as the generating function of
Betti numbers of nonsingular graded/cyclic quiver varieties.
We axiomize its properties. The axioms are purely combinatorial
statements in $\widehat{\mathscr Y}_t$, involving no geometry nor
representation theory of $\Ule$. Moreover, the axioms uniquely
characterize $\widehat\chi_{\ve,t}$, and give us an algorithm for
computation. Therefore the axioms can be considered as a definition of
$\widehat\chi_{\ve,t}$. When $\g$ is not of type $E_8$, we can
directly prove the existence of $\widehat\chi_{\ve,t}$ satisfying the
axioms without using geometry or representation theory of $\Ule$.

Two of the axioms are most important. One is the characterization of
the image of $\widehat\chi_{\ve,t}$. Another is the multiplicative
property.

The former is a modification of Frenkel-Mukhin's result~\cite{FM}.
They give a characterization of the image of $\chi_\ve$, as an analog
of the Weyl group invariance of the ordinary character homomorphism.
And they observed that the characterization gives an algorithm
computing $\chi_\ve$ at {\it l\/}--fundamental representations. This
property has no counterpart in the ordinary character homomorphism for
$\g$, and is one of the most remarkable feature of $\chi_\ve$.
We use a $t$--analog of their characterization to ``compute''
$\widehat\chi_{\ve,t}$ for {\it l\/}--fundamental representations.

A standard module $M(P)$ is a tensor product of {\it l\/}--fundamental
representations in $\operatorname{Rep}\Ule$ (see \corref{cor:product}
or \cite{VV-std}). If $\widehat\chi_{\ve,t}$ would be a ring
homomorphism, then $\widehat\chi_{\ve,t}(M(P))$ is just a product of
$\widehat\chi_{\ve,t}$ of {\it l\/}--fundamental representations.
This is {\it not\/} true under usual ring structures on $\bfR_t$ and
$\widehat{\mathscr Y}_t$. We introduce `twistings' of multiplications
on $\bfR_t$, $\widehat{\mathscr Y}_t$ so that $\widehat\chi_{\ve,t}$
is a ring homomorphism. The resulted algebras are {\it not\/}
commutative.

We can add another column to the table above by \cite{Lu:can}.
\begin{center}
\begin{tabular}{|c|}
\hline
$\mathbf U^-_q$: the $-$ part of the quantized enveloping algebra
\\
PBW basis
\\
canonical basis
\\
\hline
\end{tabular}
\end{center}
In fact, when $\mathfrak g$ is of type $A$, affine graded quiver
varieties are varieties used for the definition of the canonical base
\cite{Lu:can}.
Therefore it is more natural to relate $\bfR_t$ to the dual of
$\mathbf U^-_q$.
In this analogy, $\widehat\chi_{\ve,t}$ can be considered as an analog
of Feigin's map from $\mathbf U^-_q$ to the skew polynomial ring
(\cite{IM,Jo,Be,Re}).
We also have an analog of the monomial base. ($E((\mathbf c))$ in
\cite[7.8]{Lu:can}. See also \cite{CX,Re}.)

This article is organized as follows. In \secref{sec:qloop} we recall
results on quantum loop algebras and their finite dimensional
representations.
In \secref{sec:rings} we introduce a twisting of the multiplication on 
$\widehat{\mathscr Y}_t$.
In \secref{sec:axiom} we give axioms which $\widehat\chi_{\ve,t}$
satisfies and derive their consequences. In particular,
$\widehat\chi_{\ve,t}$ is uniquely determined from the axioms.
In \secref{sec:quiver} we introduce graded and cyclic quiver
varieties, which will be used to prove the existence of 
$\widehat\chi_{\ve,t}$ satisfying the axioms.
In \S\ref{sec:proof1}, \ref{sec:proof2}, \ref{sec:proof3} we check
that a generating function of Betti numbers of nonsingular
graded/cyclic quiver varieties satisfies the axioms.
In \secref{sec:perverse} we prove that the characterization of simple
modules mentioned above.
In \secref{sec:pm 1} we study the case $\ve = \pm 1$ in detail.
In \secref{sec:conjecture} we state a conjecture concerning finite
dimensional representations studied in the literatures~\cite{NT,HKOTY}.

In this introduction and also in the main body of this article, we put
`` '' to the word `compute'. What we actually do in this article is to
give a purely combinatorial algorithm to compute somethings. The
author wrote a computer program realizing the algorithm for computing
$\widehat\chi_{\ve,t}$ for {\it l\/}--fundamental representations when
$\g$ is of type $E$. Up to this moment (2001, April), the program
produces the answer except two {\it l\/}--fundamental representations
of $E_8$. It took 3 days for the last successful one, and the remaing
ones are inaccessible so far. In this sense, our character formula is
{\it not\/} computable in a strict sense.

The result of this article for generic $\ve$ was announced in
\cite{Na-ann}.

\subsection*{Acknowledgement}
The author would like to thank D.~Hernandez and E.~Frenkel for pointing out
mistakes in an earlier version of this paper.

\section{Quantum loop algebras}\label{sec:qloop}

\subsection{Definition}
Let $\g$ be a simple Lie algebra of type $ADE$ over $\C$. Let $I$ be
the index set of simple roots.
Let $\{\alpha_i\}_{i\in I}$, $\{ h_i\}_{i\in I}$, $\{ \Lambda_i
\}_{i\in I}$ be the sets of simple roots, simple coroots and
fundamental weights of $\g$ respectively. Let $P$ be the weight
lattice, and $P^*$ be its dual. Let $P^+$ be the semigroup of dominant 
weights.

Let $q$ be an indeterimate. For nonnegative integers $n\ge r$, define
\begin{equation*}
  [n]_q \defeq \frac{q^n - q^{-n}}{q - q^{-1}}, \quad
  [n]_q ! \defeq 
  \begin{cases}
   [n]_q [n-1]_q \cdots [2]_q [1]_q &(n > 0),\\
   1 &(n=0),
  \end{cases}
  \quad
  \begin{bmatrix}
  n \\ r
  \end{bmatrix}_q \defeq \frac{[n]_q !}{[r]_q! [n-r]_q!}.
\end{equation*}
Later we consider another indeterminate $t$. We define a $t$-binomial
coefficient
\(
\left[\begin{smallmatrix}
  n \\ r
\end{smallmatrix}\right]_t
\)
by replacing $q$ by $t$.

Let $\Ul$ be the quantum loop algebra associated with the loop algebra
$\Lg = \g\otimes \C[z,z^{-1}]$ of $\g$. It is an associative
$\Q(q)$-algebra generated by $e_{i,r}$, $f_{i,r}$ ($i\in I$,
$r\in\Z$), $q^h$ ($h\in P^*$), $h_{i,m}$ ($i\in I$, $m\in
\Z\setminus\{0\}$) with the following defining relation:
{\allowdisplaybreaks[4]
\begin{subequations}
\begin{gather*}
%
  q^0 = 1, \quad q^h q^{h'} = q^{h+h'}, \quad
  [q^h, h_{i,m}] = 0, \quad
  [h_{i,m}, h_{j,n}] = 0,
\label{eq:relHH2}\\
  q^h e_{i,r} q^{-h} = q^{\langle h,\alpha_i\rangle} e_{i,r},
  \quad
  q^h f_{i,r} q^{-h} = q^{-\langle h,\alpha_i\rangle} f_{i,r},
  \label{eq:relHE'}
\\
  ( z - q^{\pm \langle h_j,\alpha_i\rangle} w)
    \psi_i^s(z) x_j^\pm(w) =
  ( q^{\pm\langle h_j,\alpha_i\rangle} z -  w)
    x_j^\pm(w) \psi_i^s(z), \quad
    \label{eq:relHE}
\\
  \left[x_{i}^+(z), x_{j}^-(w)\right] =
  \frac{\delta_{ij}}{q  - q^{-1}}
  \left\{\delta\left(\frac{w}{z}\right)\psi^+_i(w) -
        \delta\left(\frac{z}{w}\right)\psi^-_i(z)\right\},
    \label{eq:relEF}
\\
   (z - q^{\pm 2\langle h_j,\alpha_i\rangle} w)
       x_{i}^\pm(z) x_{j}^\pm(w)
   = (q^{\pm 2\langle h_j,\alpha_i\rangle} z - w
   ) x_{j}^\pm(w) x_{i}^\pm(z),
\\
  \sum_{\sigma\in S_b}
   \sum_{p=0}^{b}(-1)^p 
   \begin{bmatrix} b \\ p\end{bmatrix}_{q}
   x_{i}^\pm(z_{\sigma(1)})\cdots x_{i}^\pm(z_{\sigma(p)})
   x_{j}^\pm(w)
   x_{i}^\pm(z_{\sigma(p+1)})\cdots x_{j}^\pm(z_{\sigma(b)}) = 0,
   \quad \text{if $i\neq j$,}
 \label{eq:relDS}
\end{gather*}
\end{subequations}
where}
$s = \pm$,
$b = 1-\langle h_i, \alpha_j\rangle$, 
and 
$S_b$ is the symmetric group of $b$ letters.
Here $\delta(z)$, $x_i^+(z)$, $x_i^-(z)$, $\psi^{\pm}_{i}(z)$ are
generating functions defined by
{\allowdisplaybreaks[4]
\begin{gather*}
   \delta(z) \defeq \sum_{r=-\infty}^\infty z^{r}, \qquad
   x_i^+(z) \defeq \sum_{r=-\infty}^\infty e_{i,r} z^{-r}, \qquad
   x_i^-(z) \defeq \sum_{r=-\infty}^\infty f_{i,r} z^{-r}, \\
   \psi^{\pm}_i(z)
  \defeq q^{\pm h_i}
   \exp\left(\pm (q-q^{-1})\sum_{m=1}^\infty h_{i,\pm m} z^{\mp m}\right).
\end{gather*}
We} also need the following generating function
\begin{equation*}
   p_i^\pm(z) \defeq 
   \exp\left(
     - \sum_{m=1}^\infty \frac{h_{i,\pm m}}{[m]_{q}} z^{\mp m}
   \right).
\end{equation*}
We have
\(
  \psi^{\pm}_i(z) = q^{\pm h_i}
  p_i^\pm(q z)/p_i^\pm(q^{-1} z).
\)

Let $e_{i,r}^{(n)} \defeq e_{i,r}^n / [n]_{q}!$,
$f_{i,r}^{(n)} \defeq f_{i,r}^n / [n]_{q}!$.
Let $\Uli$ be the $\Z[q,q^{-1}]$-subalgebra generated by
$e_{i,r}^{(n)}$, $f_{i,r}^{(n)}$ and $q^h$
for $i\in I$, $r\in \Z$, $h\in P^*$.

Let $\Uli^+$ (resp.\ $\Uli^-$) be $\Z[q,q^{-1}]$-subalgebra generated
by $e_{i,r}^{(n)}$ (resp.\ $f_{i,r}^{(n)}$) for $i\in I$, $r\in \Z$,
$n\in Z_{> 0}$.
Let $\Uli^0$ be the $\Z[q,q^{-1}]$-subalgebra generated by $q^h$, the
coefficients of $p_i^\pm(z)$ and
\begin{equation*}
   \begin{bmatrix}
      q^{h_i}; n \\ r
   \end{bmatrix}
   \defeq
   \prod_{s=1}^r \frac{q^{h_i} q^{n-s+1} - q^{-h_i} q^{-n+s-1}}
   {q^s - q^{-s}}
\end{equation*}
for all $h\in P$, $i\in I$, $n\in \Z$, $r\in \Z_{> 0}$. We have
$\Uli = \Uli^+\cdot \Uli^0\cdot \Uli^-$ (\cite[6.1]{CP-roots}).

Let $\ve$ be a nonzero complex number.
The specialization $\Uli\otimes_{\Z[q,q^{-1}]}\C$ with respect to the
homomorphism $\Z[q,q^{-1}]\ni q\mapsto \varepsilon\in\C^*$ is denoted
by $\Ule$. Set
\begin{equation*}
   \Ule^\pm \defeq \Uli^\pm\otimes_{\Z[q,q^{-1}]}\C, \qquad
   \Ule^0 \defeq \Uli^0\otimes_{\Z[q,q^{-1}]}\C.
\end{equation*}

It is known that $\Ul$ is isomorphic to a subquotient of the quantum
affine algebra $\Ua$ defined in terms of Chevalley generators $e_i$,
$f_i$, $q^h$ ($i\in I\cup \{0\}$, $h\in P^*\oplus\Z c$). (See
\cite{Dr,Be}.)
Using this identification, we define a coproduct on $\Ul$ by
\begin{equation*}\label{eq:comul}
\begin{gathered}
   \Delta q^h = q^h \otimes q^h, \quad
   \Delta e_i = e_i\otimes q^{-h_i} + 1 \otimes e_i,
\\
   \Delta f_i = f_i\otimes 1 + q^{h_i} \otimes f_i.
\end{gathered}
\end{equation*}
Note that this is different from one in \cite{Lu-book}, although there
is a simple relation between them \cite[1.4]{Kas}.
The results in \cite{Na-qaff} hold for either comultiplication
(tensor products appear in (1.2.19) and (14.1.2)).
In \cite[\S2]{Na-ann} another comultiplication was used.

It is known that the subalgebra $\Uli$ is preserved under
$\Delta$. Therefore $\Ule$ also has an induced coproduct.

For $a\in\C^*$, there is a Hopf algebra automorphism $\tau_a$ of $\Ul$ 
given by
\begin{equation*}
   \tau_a(e_{i,r}) = a^r e_{i,r}, \quad
   \tau_a(f_{i,r}) = a^r f_{i,r}, \quad
   \tau_a(h_{i,m}) = a^{m} h_{i,m}, \quad
   \tau_a(q^h) = q^h.
\end{equation*}
It preserves $\Uli\otimes_{\Z[q,q^{-1}]}\C[q,q^{-1}]$ and induces an
autormorphism of $\Ule$, which is denoted also by $\tau_a$.

We define an algebra homomorphism from ${\mathbf U}_{\ve}(\g)$ to $\Ule$ by
\begin{equation}\label{eq:subalg}
  e_i\mapsto e_{i,0},\quad f_i\mapsto f_{i,0}, \quad
  q^h \mapsto q^h,  \qquad (i\in I, h\in P^*).
\end{equation}
(See \cite[\S1.1]{Na-qaff} for the definition of ${\mathbf U}_{\ve}(\g)$.)

\subsection{Finite dimensional representation of $\Ule$}

Let $V$ be a $\Ule$-module. For $\lambda\in P$, we define
\begin{equation*}
   V_\lambda \defeq \left\{ v\in V\left|\,
   q^{h} v = \ve^{\langle h,\lambda\rangle} v,
   \begin{bmatrix}
      q^{h_i}; 0 \\ r
   \end{bmatrix} v = 
   \begin{bmatrix}
      \langle h_i, \lambda\rangle \\ r
   \end{bmatrix}_{\ve} v \right\}\right..
\end{equation*}
The module $V$ is said to be of {\it type $1$\/} if $V =
\bigoplus_\lambda V_\lambda$.
In what follows we consider only modules of type $1$.

By \eqref{eq:subalg} any $\Ule$-module $V$ can be considered as a
${\mathbf U}_{\ve}(\g)$-module. This is denoted by
$\operatorname{Res}V$. The above definition is based on the definition 
of type $1$ representation of ${\mathbf U}_{\ve}(\g)$, i.e., $V$ is of 
type $1$ if and only if $\operatorname{Res}V$ is of type $1$.

A $\Ule$-module $V$ is said to be an {\it l--highest weight module\/}
if there exists a vector $v$ such that
\(
   \Ule^+\cdot v = 0
\),
\(
   \Ule^0\cdot v \subset \C v
\)
and
\(
   V = \Ule\cdot v
\).
Such $v$ is called an {\it l--highest weight vector}.

\begin{Theorem}[\cite{CP-roots}]
  A simple {\it l\/}--highest weight module $V$ with an {\it
    l\/}--highest weight vector $v$ is finite dimensional if and only
  if there exists an $I$--tuple of polynomials $P = (P_i(u))_{i\in I}$
  with $P_i(0) = 1$ such that
\begin{gather*}
   q^h v = \ve^{\langle h, \sum_i \deg P_i \Lambda_i\rangle} v,
\quad
   \begin{bmatrix}
      q^{h_i}; 0 \\ r
   \end{bmatrix} v = 
   \begin{bmatrix}
      \deg P_i \\ r
   \end{bmatrix}_{\ve} v,
\\
   p_i^+(z) v = P_i(1/z) v, \quad
   p_i^-(z) v = c_{P_i}^{-1} z^{\deg P_i} P_i(1/z) v,
\end{gather*}
where $c_{P_i}$ is the top term of $P_i$, i.e., the coefficient of
$u^{\deg P_i}$ in $P_i$.
\end{Theorem}

The $I$-tuple of polynomials $P$ is called the {\it l--highest
  weight}, or the Drinfeld polynomial of $V$.
We denote the above module $V$ by $L(P)$ since it is 
determined by $P$.

For $i\in I$ and $a\in\C^*$, the simple module $L(P)$ with
\begin{equation*}
   P_i(u) = 1 - au, \qquad
   P_j(u) = 1 \quad\text{if $j\neq i$}
\end{equation*}
is called an {\it l--fundamental representation\/} and denoted by
$L(\Lambda_i)_a$.

Let $V$ be a finite dimensional $\Ule$-module with the weight space
decomposition $V = \bigoplus V_\lambda$. Since the commutative
subalgebra $\Ule^0$ preserves each $V_\lambda$, we can further
decompose $V$ into a sum of generalized simultaneous eigenspaces of
$\Ule^0$.

\begin{Theorem}[{\protect \cite[Proposition 1]{FR}, \cite[Lemma
3.1]{FM2}, \cite[13.4.5]{Na-qaff}}]
Simultaneous eigenvalues of $\Ule^0$ have the following forms:
\begin{gather*}
   \ve^{\langle h, \deg Q^1_i - \deg Q^2_i\rangle}\quad \text{for $q^h$},
\qquad
   \begin{bmatrix}
     \deg Q^1_i - \deg Q^2_i \\ r
   \end{bmatrix}_\ve \quad \text{for $\begin{bmatrix}
      q^{h_i}; 0 \\ r
\end{bmatrix}$},
\\
   \frac{Q^1_i(1/z)}{Q^2_i(1/z)}\quad \text{for $p_i^+(z)$},
\qquad
   \frac{c_{Q^1_i}^{-1}z^{\deg Q^1_i} Q^1_i(1/z)}
   {c_{Q^2_i}^{-1}z^{\deg Q^2_i}Q^2_i(1/z)}\quad \text{for $p_i^-(z)$},
\end{gather*}
where $Q^1_i$, $Q^2_i$ are polynomials with $Q^1_i(0) = Q^2_i(0) = 1$
and $c_{Q^1_i}$, $c_{Q^2_i}$ are as above.
\end{Theorem}

We simply write the $I$-tuple of rational functions
$(Q^1_i(u)/Q^2_i(u))$ by $Q$.
A generalized simultaneous eigenspace is called an {\it l--weight
  space}. The corresponding $I$-tuple of rational functions is called
an {\it l--weight}. We denote the {\it l\/}--weight space by $V_Q$.

The $q$--character, or $\ve$--character \cite{FR,FM2} of a finite
dimensional $\Ule$-module $V$ is defined by
\begin{equation*}
   \chi_\ve(V) = \sum_Q \dim V_Q\; e^Q.
\end{equation*}
The precise definition of $e^Q$ will be explained in the next section.

\subsection{Standard modules}

We will use another family of finite dimensional {\it l\/}--highest weight
modules, called standard modules.

Let $\bw\in P^+$ be a dominant weight. Let $w_i = \langle
h_i,\bw\rangle\in\Z_{\ge 0}$.
Let $G_{\bw} = \prod_{i\in I} \GL(w_i,\C)$.
Its representation ring
$R(G_{\bw})$ is the invariant part of the Laurant polynomial ring:
\begin{equation*}
   R(G_{\bw})
   = \Z[x_{1,1}^\pm,\dots, x_{1,w_1}^\pm]^{\mathfrak S_{w_1}}
   \otimes \Z[x_{2,1}^\pm,\dots,
     x_{2,w_2}^\pm]^{\mathfrak S_{w_2}}
   \otimes\cdots\otimes
     \Z[x_{n,1}^\pm,\dots, x_{n,w_n}^\pm]^{\mathfrak S_{w_n}},
\end{equation*}
where we put a numbering $1,\dots,n$ to $I$.
In \cite{Na-qaff}, we constructed a $\Uli\otimes_{\Z} R(G_\bw)$-module
$M(\bw)$ such that it is free of finite rank over
$R(G_\bw)\otimes\Z[q,q^{-1}]$ and has a vector $[0]_{\bw}$ satisfying
\begin{subequations}
\begin{gather*}
   e_{i,r}[0]_\bw = 0\quad\text{for any $i\in I$, $r\in \Z$},
\\
   M(\bw) = \left(\Uli^-\otimes_{\Z} R(G_\bw)\right)[0]_\bw,
   \label{eq:span}
\\
   q^h [0]_\bw = q^{\langle h,\bw\rangle} [0]_\bw,
\\
   p_i^+(z)[0]_\bw
   = \prod_{p=1}^{w_i} \left(1-\frac{x_{i,p}}{z}\right)[0]_\bw,
\\
   p_i^-(z)[0]_\bw
   = \prod_{p=1}^{w_i} \left(1-\frac{z}{x_{i,p}}\right)[0]_\bw,
\end{gather*}
\end{subequations}
If an $I$-tuple of monic polynomials $P(u) =
(P_i(u))_{i\in I}$ with $\deg P_i = w_i$ is given,
then we define a {\it standard module\/} by the specialization
\begin{equation*}
   M(P) = M(\bw)\otimes_{R(G_\bw)[q,q^{-1}]} \C,
\end{equation*}
where the algebra homomorphism $R(G_\bw)[q,q^{-1}]\to \C$ sends
$q$ to $\ve$ and $x_{i,1},\dots, x_{i,w_k}$ to roots of $P_i$.
The simple module $L(P)$ is the simple quotient of $M(P)$.

The original definition of the universal standard module
\cite{Na-qaff} is geometric. However, it is not difficult to give an
algebraic characterization.  Let $M(\Lambda_i)$ be the universal
standard module for the dominant weight $\Lambda_i$. It is a
$\Uli[x,x^{-1}]$-module. Let $W(\Lambda_i) =
M(\Lambda_i)/(x-1)M(\Lambda_i)$. Then we have
\begin{Theorem}[{\protect\cite[1.22]{Na-tensor}}]
Put a numbering $1,\dots,n$ on $I$.  Let $w_i = \langle
h_i,\bw\rangle$.  The universal standard module $M(\bw)$ is the
$\Uli\otimes_\Z R(G_\lambda)$-submodule of
\[
   W({\Lambda_1})^{\otimes w_1}
   \otimes \cdots\otimes
   W({\Lambda_n})^{\otimes w_n}
   \otimes
   \Z[q,q^{-1},x_{1,1}^\pm,\dots, x_{1,w_1}^\pm,
    \cdots,
    x_{n,1}^\pm,\dots, x_{n,w_n}^\pm]
\]
\rom(the tensor product is over $\Z[q,q^{-1}]$\rom)
generated by
\(
   \bigotimes_{i\in I} [0]_{\Lambda_i}^{\otimes \lambda_i}
\).
\rom(The result holds for the tensor product of any order.\rom)
\end{Theorem}

It is not difficult to show that $W(\Lambda_i)$ is isomorphic to a
module studied by Kashiwara~\cite{Kas2} ($V(\lambda)$ in his
notation). Since his construction is algebraic, the standard module
$M(\bw)$ has an algebraic construction.

We also prove that $M(P^1P^2)$ is equal to $M(P^1)\otimes M(P^2)$ in
the representation ring $\operatorname{Rep}\Ule$ later.  (See
\corref{cor:product}.)  Here the $I$-tuple of polynomials $(P_i
Q_i)_i$ for $P = (P_i)_i$, $Q = (Q_i)_i$ is denoted by $PQ$ for
brevity.

\section{A modified multiplication on $\widehat{\mathscr Y}_t$}
\label{sec:rings} 

We use the following polynomial rings in this article:
{\allowdisplaybreaks[4]
\begin{equation*}
\begin{split}
  & \widehat{\mathscr Y}_t \defeq
  \Z[t, t^{-1}, V_{i,a}, W_{i,a}]_{i\in I,a\in\C^*},
\\
  & \mathscr Y_t \defeq 
   \Z[t,t^{-1},Y_{i,a}, Y_{i,a}^{-1}]_{i\in I, a\in\C^*},
\\
  & \mathscr Y \defeq
   \Z[Y_{i,a}, Y_{i,a}^{-1}]_{i\in I, a\in\C^*},
\\
  & \overline{\mathscr Y} \defeq
  Z[y_i,y_i^{-1}]_{i\in I}.   
\end{split}
\end{equation*}
}

We consider $\widehat{\mathscr Y}_t$ as a polynomial ring in
infinitely many variables $V_{i,a}$, $W_{i,a}$ with coefficients in
$\Z[t,t^{-1}]$. So a {\it monomial\/} means a monomial only in
$V_{i,a}$, $W_{i,a}$, containing no $t$, $t^{-1}$. The same
convention applies also to $\mathscr Y_t$.

For a monomial $m\in\widehat{\mathscr Y}_t$, let $w_{i,a}(m)$,
$v_{i,a}(m)\in\Z_{\ge 0}$ be the degrees in $V_{i,a}$, $W_{i,a}$,
i.e.,
\begin{equation*}
   m = \prod_{i,a} V_{i,a}^{v_{i,a}(m)}\, W_{i,a}^{w_{i,a}(m)}.
\end{equation*}
We also define
\begin{equation*}
   u_{i,a}(m) \defeq
   w_{i,a}(m) - v_{i,a\ve^{-1}}(m) - v_{i,a\ve}(m)
 + \sum_{j:C_{ji} = -1} v_{j,a}(m).
\end{equation*}

When $\ve$ is not a root of unity, we define $(\tilde
u_{i,a}(m))_{i\in I, a\in\C^*}$ for a monomial $m$ in
$\widehat{\mathscr Y}_t$, as the solution of
\begin{equation*}
   u_{i,a}(m) = \tilde u_{i,a\ve^{-1}}(m) + \tilde u_{i,a\ve}(m)
                  - \sum_{j: a_{ij} = -1} \tilde u_{j,a}(m).
\end{equation*}
To solve the system, we may assume that $u_{i,a}(m)=0$ unless $a$ is a
power of $q$. Then the above is a recursive system, since $q$ is not a
root of unity. So it has a unique solution such that $\tilde
u_{i,q^s}(m) = 0$ for sufficiently small $s$. Note that $\tilde
u_{i,a}(m)$ is nonzero for possibly infinitely many $a$'s, although
$u_{i,a}(m)$ is not.

If $m^1$, $m^2$ are monomials, we set
\begin{equation}\label{eq:def_d}
  \begin{split}
   d(m^1, m^2) & \defeq 
   \sum_{i,a} \left( v_{i,a\ve}(m^1) u_{i,a}(m^2)
   + w_{i,a\ve}(m^1) v_{i,a}(m^2)\right)
\\
   &= \sum_{i,a} \left( u_{i,a}(m^1) v_{i,a\ve^{-1}}(m^2) 
   + v_{i,a}(m^1) w_{i,a\ve^{-1}}(m^2) \right).
  \end{split}
\end{equation}
From the definition, $d(\ ,\ )$ satisfies
\begin{equation}\label{eq:shift}
   d(m^1 m^2, m^3) = d(m^1, m^3) + d(m^2, m^3), \qquad
   d(m^1, m^2 m^3) = d(m^1, m^2) + d(m^1, m^3).
\end{equation}

When $\ve$ is not a root of unity, we also define
\begin{equation*}
   \tilde d(m^1, m^2) \defeq 
   - \sum_{i,a} u_{i,a}(m^1) \tilde u_{i,a\ve^{-1}}(m^2)
.
\end{equation*}
Since $u_{i,a}(m^2) = 0$ except for finitely many $a$'s, this is
well-defined. Moreover, we have
\begin{equation*}
   \tilde d(m^1, m^2) = d(m^1, m^2) + \tilde d_W(m^1, m^2),
\end{equation*}
where $\tilde d_W$ is defined as $\tilde d$ by replacing $u_{i,a}$ by
$w_{i,a}$. Here we have used $\tilde u_{i,a}(m) = \tilde w_{i,a}(m) -
v_{i,a}(m)$.

We define an involution
$\setbox5=\hbox{A}\overline{\rule{0mm}{\ht5}\hspace*{\wd5}}$ on
$\widehat{\mathscr Y}_t$ by 
\begin{equation}
\label{eq:bar}
   \overline{t} = t^{-1}, \quad
   \overline{m} = t^{2d(m,m)} m,
\end{equation}
where $m$ is a monomial in $V_{i,a}$, $W_{i,a}$.
We define an involution
$\setbox5=\hbox{A}\overline{\rule{0mm}{\ht5}\hspace*{\wd5}}$ on
${\mathscr Y}_t$ by $\overline{t} = t^{-1}$,
$\overline{Y_{i,a}^\pm} = Y_{i,a}^\pm$.

We define a new multiplication $\ast$ on $\widehat{\mathscr Y}_t$ by
\begin{equation*}
   m^1 \ast m^2 \defeq t^{2d(m^1,m^2)} m^1 m^2,
\end{equation*}
where $m^1$, $m^2$ are monomials and $m^1 m^2$ is the usual
multiplication of $m^1$ and $m^2$. By \eqref{eq:shift} it is
associative.
({\bf NB}: The multiplication in \cite{Na-ann} was $m^1 \ast m^2
\defeq t^{2d(m^2,m^1)} m^1 m^2$. This is because the coproduct is
changed.)

From definition we have
\begin{equation}\label{eq:antihom}
  \overline{m^1\ast m^2} = \overline{m^2}\ast\overline{m^1}.
\end{equation}

Let us give an example which will be important later. Suppose that
$m$ is a monomial with $u_{i,a}(m) = 1$, $u_{i,b}(m) = 0$ for
$b\neq a$. Then
\begin{equation}\label{eq:example}
\begin{split}
   \left[m(1 + V_{i,a\ve})\right]^{\ast n}
   & \defeq \underbrace{
   m(1 + V_{i,a\ve})\ast \cdots
   \ast m(1 + V_{i,a\ve})}_{\text{$n$ times}}
\\
   & = m^n \sum_{r=0}^n 
    t^{r(n-r)}
   \begin{bmatrix}
     n \\ r
   \end{bmatrix}_t V_{i,a\ve}^r.
\end{split}
\end{equation}
When $\ve$ is not a root of unity, we define another multiplication
$\tilde\ast$ by
\begin{equation*}
   m^1 \tilde\ast\, m^2 
   \defeq t^{\tilde d(m^1,m^2)-\tilde d(m^2,m^1)} m^1 m^2.
\end{equation*}

We define a $\Z[t,t^{-1}]$-linear homomorphism
\(
   \widehat\Pi\colon\widehat{\mathscr Y}_t\to \mathscr Y_t
\)
by
\begin{equation}\label{eq:YY}
  m = \prod_{i,a} V_{i,a}^{v_{i,a}(m)}\, W_{i,a}^{w_{i,a}(m)}
    \longmapsto t^{-d(m,m)} \prod_{i,a} Y_{i,a}^{u_{i,a}(m)}.
\end{equation}
This is not a ring homomorphism with respect to either the ordinary
multiplication or $\ast$. However, when $\ve$ is not a root of unity,
we can define a new multiplication on $\mathscr Y_t$ so that the above
is a ring homomorphism with respect to this multiplication and $\tilde\ast$.
It is because $\ve(m^1, m^2)$ involves only $u_{i,a}(m^1)$, $u_{i,a}(m^2)$.
We denote also by $\tilde\ast$ the new multiplication on $\mathscr
Y_t$. We have
\begin{equation}\label{eq:use}
\begin{gathered}
\widehat\Pi(m^1 \ast m^2)
   = t^{\tilde d_W(m^1,m^2) - \tilde d_W(m^2, m^1)}
   \widehat\Pi(m^1) \tilde\ast \widehat\Pi(m^2),
\\
   \widehat\Pi\circ
     \setbox5=\hbox{A}\overline{\rule{0mm}{\ht5}\hspace*{\wd5}}
   = \setbox5=\hbox{A}\overline{\rule{0mm}{\ht5}\hspace*{\wd5}}
     \circ\widehat\Pi.
\end{gathered}
\end{equation}

Further we define homomorphisms $\Pi_t\colon\mathscr Y_t\to\mathscr Y$,
$\overline{\Pi}\colon\mathscr Y \to \Z[y_i^\pm]$ by
\begin{equation*}
   \Pi_t\colon \mathscr Y_t \ni 
   \begin{aligned}[c]
    t & \longmapsto 1
   \\
    Y_{i,a} & \longmapsto Y_{i,a}
   \end{aligned}
   \in \mathscr Y, \qquad
   \overline\Pi\colon 
   \mathscr Y\ni Y_{i,a} \longmapsto y_i\in \Z[y_i, y_i^{-1}]_{i\in I}.
\end{equation*}
The composition $\widehat{\mathscr Y}_t\to\mathscr Y$ or
$\widehat{\mathscr Y}_t\to\Z[y_i^\pm]$ is a ring homomorphism with
respect to both the usual multiplication and $\ast$.

\begin{Definition}
A monomial $m\in\widehat{\mathscr Y}_t$ is said {\it $i$--dominant\/} if
$u_{i,a}(m)\ge 0$ for any $i\in I$.
A monomial $m\in\widehat{\mathscr Y}_t$ is said {\it l--dominant\/} if
it is $i$--dominant for all $i\in I$, i.e., $\widehat\Pi(m)$ contains
only nonnegative powers of $Y_{i,a}$.
Similarly a monomial $m\in\mathscr Y$ is said {\it l--dominant\/} if
it contains only nonnegative powers of $Y_{i,a}$.
Note that a monomial $m\in\Z[y_i, y_i^{-1}]_{i\in I}$ contains only
nonnegative powers of $y_{i}$ if and only if it is dominant as a
weight of $\g$.
\end{Definition}

Let
\[
   m = \prod_{i,a} Y_{i,a}^{u_{i,a}}
\]
be a monomial in $\mathscr Y$ with $u_{i,a}\in\Z$. We associate to $m$
an $I$-tuple of rational functions $Q = (Q_i)$ by
\begin{equation*}
   Q_i(u) = \prod_a \left(1 - au\right)^{u_{i,a}}.
\end{equation*}
Conversely an $I$-tuple of rational functions $Q = (Q_i)$ as above
determines a monomial in $\mathscr Y$. We denote it by $e^Q$.
This is the $e^Q$ mentioned in the previous section. Note that $e^Q$
is {\it l\/}--dominant if and only if $Q$ is an $I$-tuple of
polynomials.

We also use a similar identification between an $I$-tuple of
polynomials $P = (P_i)$ and a monomial $m$ in $W_{i,a}$ ($i\in I$,
$a\in \C^*$):
\begin{equation*}
   m = \prod_{i,a} W_{i,a}^{w_{i,a}} \longleftrightarrow
   P = (P_i);\quad P_i(u) = \prod_a (1 - au)^{w_{i,a}}.
\end{equation*}
We denote $m$ also by $e^P$, hoping that it makes no confusion.

\begin{Definition}
Let $m$, $m'$ be monomials in $\widehat{\mathscr Y}_t$.
We say $m\le m'$ if $m/m'$ is a monomial in $V_{i,a}$ ($i\in I$,
$a\in\C^*$). We say $m < m'$ if $m \le m'$ and $m\neq m'$. It defines
a partial order among monomials in $\widehat{\mathscr Y}_t$.
Similarly for monomials $m$, $m'$ in $\mathscr Y$, we say 
$m\le m'$ if $m/m'$ is a monomial in 
\(
   {\widehat\Pi}(V_{i,a})
\)
($i\in I$, $a\in\C^*$).
For two $I$-tuple of rational functions $Q$, $Q'$, we say $Q\le Q'$ if
$e^Q \le e^{Q'}$.
Finally for monomials $m$, $m'$ in $\Z[y_i,y_i^{-1}]_{i\in I}$, we say 
$m\le m'$ if $m/m'$ is a monomial in 
\(
   {\overline{\Pi}\circ\Pi_t\circ\widehat\Pi}(V_{i,a})
\)
($i\in I$, $a\in\C^*$). But this is nothing but the usual order on
weights.
\end{Definition}

\section{A $t$--analog of the $q$--character: Axioms}\label{sec:axiom}

A main tool in this article is a $t$--analog of the $q$--character:
\begin{equation*}
   \widehat\chi_{\ve,t}\colon \bfR_t
   = \operatorname{Rep}\Ule\otimes_\Z \Z[t,t^{-1}]
   \to \widehat{\mathscr Y}_t.
\end{equation*}
For the definition we need geometric constructions of standard
modules, so we will postpone it to \secref{sec:quiver}.
In this section, we explain properties of $\widehat\chi_{\ve,t}$ as
axioms. Then we show that these axioms uniquely characterize
$\widehat\chi_{\ve,t}$, and in fact, give us an algorithm for
``computation''.
Thus we may consider the axioms as the definition of
$\widehat\chi_{\ve,t}$.

Our first axiom is the highest weight property:
\begin{Axiom}
The $\widehat\chi_{\ve,t}$ of a standard module $M(P)$ has a form
\begin{equation*}
   \widehat\chi_{\ve,t}(M(P)) = e^P + \sum a_m(t) m,
\end{equation*}
where each monomial $m$ satisfies $m < e^P$.
\end{Axiom}

Composing maps $\widehat{\mathscr Y}_t\to\mathscr Y_t$,
$\mathscr Y_t\to\mathscr Y$, $\mathscr Y\to \Z[y_i^\pm]$ in
\secref{sec:rings}, we define maps
\begin{gather*}
    \chi_{\ve,t}=\widehat\Pi\circ\widehat\chi_{\ve,t}
      \colon \bfR_t \to \mathscr Y_t,
\\
    \chi_\ve=\Pi_t\circ\chi_{\ve,t}
     \colon\operatorname{Rep}\Ule \to \mathscr Y,
\qquad
    \chi=\overline{\Pi}\circ\chi_\ve
     \colon\operatorname{Rep}\Ule \to \Z[y_i, y_i^{-1}]_{i\in I}.
\end{gather*}
The $\widehat\chi_{\ve,t}$ is a homomorphism of
$\Z[t,t^{-1}]$-modules, not of rings.

Frenkel-Mukhin \cite[5.1, 5.2]{FM} proved that the image of $\chi_\ve$ is
equal to
\begin{equation*}
  \bigcap_{i\in I}
  \left( \Z[Y_{j,a}^\pm]_{j:j\neq i, a\in\C^*}
    \otimes\Z[Y_{i,b}(1+V_{i,b\ve})]_{b\in\C^*}\right).
\end{equation*}
We define its $t$--analog, replacing $(1+V_{i,b\ve})^n$ by
\begin{equation*}
  \sum_{r=0}^n t^{r(n-r)}\begin{bmatrix} n \\ r
  \end{bmatrix}_t V_{i,b\ve}^{r}.
\end{equation*}
More precisely, for each $i\in I$, let $\widehat{\mathscr K}_{t,i}$ be
the $\Z[t,t^{-1}]$-linear subspace of $\widehat{\mathscr Y}_t$
generated by elements
\begin{equation}\label{eq:form}
   E_i(m) \defeq 
    m\, \prod_a
     \sum_{r_a=0}^{u_{i,a}(m)}
     t^{r_a(u_{i,a}(m)-r_a)}
     \begin{bmatrix}
       u_{i,a}(m) \\ r_a
     \end{bmatrix}_t V_{i,a\ve}^{r_a},
\end{equation}
where $m$ is an $i$--dominant monomial, i.e., $u_{i,a}(m) \ge 0$ for
all $a\in\C^*$.
Let
\[
   \widehat{\mathscr K}_t \defeq
    \bigcap_i \widehat{\mathscr K}_{t,i}, \qquad
   \mathscr K_t \defeq \widehat\Pi(\widehat{\mathscr K}_t)
   \subset \mathscr Y_t.
\]

\begin{Axiom}
%
The image of $\widehat\chi_{\ve,t}$ is contained in $\widehat{\mathscr
K}_t$.
\end{Axiom}

Next axiom is about the multiplicative property of
$\widehat\chi_{\ve,t}$. As explained in the introduction, it is not
multiplicative under the usual product structure on $\bfR_t$.

\begin{Axiom}
Suppose that two $I$-tuples of polynomials $P^1 = (P^1_i)$, $P^2 =
(P^2_i)$ satisfy the following condition:
\begin{equation}
\label{eq:Z}
\begin{minipage}[m]{0.75\textwidth}
\noindent   
   $a/b\notin\{ \ve^n \mid n\in\Z, n \ge 2\}$ for any
   pair $a$, $b$ with $P^1_i(1/a) = 0$, $P^2_j(1/b) =
   0$ \textup($i,j\in I$\textup).
\end{minipage}
\end{equation}
Then we have
\begin{equation*}
  \widehat\chi_{\ve,t}(M(P^1P^2)) = \widehat\chi_{\ve,t}(M(P^1))\ast
  \widehat\chi_{\ve,t}(M(P^2)).
\end{equation*}
\end{Axiom}

We have the following special case
\begin{equation*}
   \widehat\chi_{\ve,t}(M(P^1P^2)) 
    = \widehat\chi_{\ve,t}(M(P^1))\widehat\chi_{\ve,t}(M(P^2))   
\end{equation*}
under the stronger condition $a/b\notin \ve^\Z$ by the
definition of $\ast$.

The last axiom is about a specialization at a root of unity.
Suppose that $\ve$ is a primitive $s$-th root of unity.
We choose and fix $q$, which is not root of unity. The axiom will say
that $\widehat\chi_{\ve,t}(M(P))$ can be written in terms of
$\widehat\chi_{q,t}(M(P_q))$ for some $P_q$.

By Axiom~3, more precisely the sentence following Axiom~3, we may
assume that inverse of roots of $P_i(u) = 0$ ($i\in I$) is contained
in $a\ve^\Z$ for some $a\in\C^*$. Therefore
\begin{equation*}
   P_i(u) = \prod_{n=0}^{s-1} (1 - a \ve^n u)^{N_{i,n}}.
\end{equation*}
with $N_{i,n}\in\Z_{\ge 0}$. We define $P_q = ((P_q)_i)$ by
\begin{equation*}
   (P_q)_i(u) = \prod_{n=0}^{s-1} (1 - a q^n u)^{N_{i,n}}.
\end{equation*}
We set $N_{i,n} = 0$ if $n\notin \{0,\dots,s-1\}$.

Let 
\[
   \widehat\chi_{q,t}(M(P_q)) = \sum a_m(t) m.
\]
By previous axioms, each $m$ is written as
\begin{equation}\label{eq:m}
   m = e^{P_q} \prod_{i\in I, n\in \Z} V_{i,a q^n}^{M_{i,n}}
     = \prod_{i\in I, n\in \Z} W_{i,a q^n}^{N_{i,n}}
                               V_{i,a q^n}^{M_{i,n}}
\end{equation}
with $M_{i,n}\in \Z_{\ge 0}$.
By previous axioms $M_{i,n}$ is independent of $q$ (cf.\ 
\thmref{thm:cons}(4)).
We define monomials $m|_{q = \ve}$, $m[k]$ by
\begin{equation}\label{eq:m_shift}
\begin{split}
   m|_{q = \ve} & \defeq
   \prod_{i\in I, n\in \Z} W_{i,a\ve^n}^{N_{i,n}}
     V_{i,a\ve^n}^{M_{i,n}},
\\
   m[k] & \defeq 
   \prod_{i\in I, n\in \Z} W_{i,a q^n}^{N_{i,n+k}}
   V_{i,a q^n}^{M_{i,n+k}}.
\end{split}
\end{equation}
Note that $m|_{q=\ve} = m[k]|_{q=\ve}$ if $k\equiv 0\mod s$.
We define
\begin{equation*}
   D^-(m) \defeq \sum_{k < 0} d_q(m,m[ks]),
\end{equation*}
where we define $d_q$ as $d$ in \eqref{eq:def_d} replacing $\ve$ by
$q$.

\begin{Axiom}
We have
\begin{equation*}
   \widehat\chi_{\ve,t}(M(P))
   = \sum t^{2D^-(m)}\, a_m(t)\, m|_{q = \ve}.
\end{equation*}
\end{Axiom}

We can consider similar axioms for 
\(
   \chi_\ve = \Pi_t\circ\widehat\Pi\circ\widehat\chi_{\ve,t}
\).
Axioms~3,4 are simplified when $t=1$. Axiom~3 is
\(
   \chi_\ve(M(P^1P^2)) = \chi_\ve(M(P^1))\chi_\ve(M(P^2))
\).
Axiom~4 says
\(
   \chi_\ve(M(P)) = \chi_q(M(P))|_{q=\ve}
\).
The original $\chi_\ve$ defined in \cite{FR,FM2} satisfies those
axioms: Axioms~1,2 were proved in \cite[Theorem~4.1, Theorem~5.1]{FM}.
Axiom~3 was proved in \cite[Lemma~3]{FR}. Axiom~4 was proved in
\cite[Theorem~3.2]{FM2}.

Let us give few consequences of the axioms.
\begin{Theorem}\label{thm:cons}
%
\textup{(1)} The map $\chi_{\ve,t}$ \textup(and hence
also $\widehat\chi_{\ve,t}$\textup) is injective.
The image of $\chi_{\ve,t}$ is equal to ${\mathscr K}_t$.

\textup{(2)} Suppose that a $\Ule$-module $M$ has the following
property: $\widehat\chi_{\ve,t}(M)$ contains only one {\it
l\/}--dominant monomial $m_0$.
Then $\widehat\chi_{\ve,t}(M)$ is uniquely determined from $m_0$ and
the condition $\widehat\chi_{\ve,t}(M)\in\widehat{\mathscr K}_t$.

\textup{(3)} Let $m$ be an {\it l\/}--dominant monomial in $\mathscr
Y_t$, considered as an element of the dual of $\bfR_t$ by taking the
coefficient of $\chi_{\ve,t}$ at $m$. Then $\{ m \mid \text{$m$ is
  {\it l\/}--dominant}\}$ is a base of the dual of $\bfR_t$.

\textup{(4)} The $\widehat\chi_{\ve,t}$ is unique, if it exists.

\textup{(5)} $\widehat\chi_{\ve,t}(\tau_a^*(V))$ is obtained from
$\widehat\chi_{\ve,t}(V)$ by replacing $W_{i,b}$, $V_{i,b}$ by
$W_{i,ab}$, $V_{i,ab}$.

\textup{(6)} The coefficients of a monomial $m$ in
$\widehat\chi_{\ve,t}(M(P))$ is a polynomial in $t^2$.  \textup(In
fact, it will become clear that it is a polynomial in $t^2$ with
nonnegative coefficients.\textup)
\end{Theorem}

\begin{proof}
These are essentially proved in \cite{FR,FM}. So our proof is sketchy.

(1) Since $\chi_{\ve,t}(M(P))$ equals $\widehat\Pi(e^P)$ plus the sum
of lower monomials, the first assertion follows by induction on $<$.
The second assertion follows from the argument in \cite[5.6]{FM},
where we use the standard module $M(P)$ instead of simple modules.

(2) Let $m$ be a monomial appearing in $\widehat\chi_{\ve,t}(M)$,
which is not $m_0$. It is not {\it l\/}--dominant by the
assumption. By Axiom~2, $m$ appears in $E_i(m')$ for some
monomial $m'$ appearing in $\widehat\chi_{\ve,t}(M)$. In particular,
we have $m < m'$. Repeating the argument for $m'$, we have $m < m_0$.

The coefficient of $m$ in $\widehat\chi_{\ve,t}(M)$ is equal to the
sum of coefficients of $m$ in $E_i(m')$ for all possible $m'$'s. ($i$
is fixed.) Again by induction on $<$, we can determine the coefficient
inductively.

(3) By Axiom~1, the transition matrix between $\{ M(P)\}$ and the dual
base of $\{ m\}$ above is uppertriangular with diagonal entries $1$.

(4) By Axiom~4, we may assume that $\ve$ is not a root of unity.
Consider the case $P_i(u) = 1 - au$, $P_j(u) = 1$ for $j\neq i$ for
some $i$. By \cite[Corollary 4.5]{FM}, Axiom~1 implies that the
$\widehat\chi_{\ve,t}(M(P))$ for $P$ does not contains {\it
l\/}--dominant terms other than $e^P$.
(See \propref{prop:nodom} below for a geometric proof.)
In particular, $\widehat\chi_{\ve,t}(M(P))$ is uniquely determined by
above (2) in this case.
We use Axiom~3 to ``calculate'' $\widehat\chi_{\ve,t}(M(P))$ for
arbitrary $P$ as follows.
We order inverses of roots (counted with multiplicities) of $P_i(u) =
0$ ($i\in I$) as $a_1$, $a_2$, \dots, so that
$a_p/a_q\neq \ve^n$ for $n\ge 2$ if $p < q$. This is
possible since $\ve$ is not a root of unity. For each $a_p$, we
define a Drinfeld polynomial $Q^p$ by
\begin{equation*}
     Q^p_{i_p}(u) = (1 - a_p u), \qquad
     Q^p_{j}(u) = 1 \quad (j\neq i_p),
\end{equation*}
if $1/a_p$ is a root of $P_{i_p}(u) = 0$. Therefore we have $P_i =
\prod_p Q^p_i$. By our choice, we have
\begin{equation*}
    \widehat\chi_{\ve,t}(M(P))
    = \widehat\chi_{\ve,t}(M(Q^1))\ast \widehat\chi_{\ve,t}(M(Q^2))
    \ast\cdots
\end{equation*}
by Axiom~3. Each $\widehat\chi_{\ve,t}(M(Q^p))$ is uniquely determined
by the above discussion. Therefore $\widehat\chi_{\ve,t}(M(P))$ is
also uniquely determined.

(5) It is enough to check the case $V = M(P)$. In this case,
$\tau_a^*(M(P))$ is the standard module with Drinfeld polynomial
$P(au)$. The assertion follows from the axioms.

(6) This also follows from the axioms. By Axiom~4, we may assume $\ve$
is a root of unity. By Axiom~3, we may assume $M(P)$ is an {\it
  l\/}--fundamental representation. In this case, the assertion
follows from Axiom~2, since 
\(
  t^{r(n-r)}\left[
  \begin{smallmatrix}
    n \\ r 
  \end{smallmatrix}\right]_t
\)
is a polynomial in $t^2$.
\end{proof}

In \cite[\S5.5]{FM}, Frenkel-Mukhin gave an explict combinatorial
algorithm to ``compute'' $\widehat\chi_{\ve,t}(M)$ for $M$ as in (2).
We will give a geometric interpretation of their algorithm in
\secref{sec:proof1}.

By the uniqueness, we get
\begin{Corollary}\label{cor:coincide}
The $\chi_{\ve}$ coincides with the $\ve$-character
defined in \cite{FR,FM2}.
\end{Corollary}


By \cite[Theorem~3]{FR}, $\chi$ is the ordinary character of the
restriction of a $\Ule$-module to a $\Ue$-module.

As promised, we prove
\begin{Corollary}\label{cor:product}
In the representation ring $\operatorname{Rep}\Ule$, we have
\begin{equation*}
   M(P^1 P^2) = M(P^1)\otimes M(P^2)
\end{equation*}
for any $I$-tuples of polynomials $P^1$, $P^2$.
\end{Corollary}

\begin{proof}
Since $\chi_\ve$ is injective, it is enough to show that 
$\chi_\ve(M(P^1 P^2)) = \chi_\ve(M(P^1))\chi_\ve(M(P^2))$.

In fact, it is easy to prove this equality directly from the geometric
defintion in \eqref{eq:geomdef}. However, we prove it only from Axioms.

By Axiom~4, we may assume $\ve$ is not a root of unity. We order
inverses of roots (counted with multiplicities) of $P^1_iP^2_i(u) = 0$
($i\in I$) as in the proof of \thmref{thm:cons}(4). Then we have
\begin{equation*}
    \chi_\ve(M(P^1P^2)) = \prod_p \chi_\ve(M(Q^p))
\end{equation*}
by Axiom~3. The product can be taken in any order, since
$\operatorname{Rep}\Ule$ is commutative. Each $a_p$ is either the
inverse of a root of $P^1_i(u) = 0$ or $P^2_i(u)=0$. We divide
$a_p$'s into two sets accordingly. Then the products of
$\chi_\ve(M(Q^a))$ over groups are equal to $\chi_\ve(M(P^1))$ and
$\chi_\ve(M(P^2))$ again by Axiom~3. Therefore we get the assertion.
\end{proof}

We also give another consequence of the axioms.

\begin{Theorem}\label{thm:bar}
%
  The $\widehat{\mathscr K}_t$ is invariant under the multiplication
  $\ast$ and the involution
  $\setbox5=\hbox{A}\overline{\rule{0mm}{\ht5}\hspace*{\wd5}}$ on
  $\widehat{\mathscr Y}_t$.
  Moreover, $\bfR_t$ has an involution induced from one on
  $\widehat{\mathscr Y}_t$. When $\ve$ is not a root of unity, it also
  has a multiplication induced from that on $\mathscr Y_t$.
\end{Theorem}

The following proof is elementary, but less conceputal. We will give
another geometric proof in \secref{sec:proof2}.

\begin{Remark}
The multiplication on $\bfR_t$ in an earlier version was not
associative, although it works for the computation of tensor product
decompositions of two simple modules.
A modification of the multiplication here was inspired by a paper of
Varagnolo-Vasserot \cite{VV2}.
\end{Remark}

\begin{proof}
For simplicity, we assume that $\ve$ is not a root of unity. The proof 
for the case when $\ve$ is a root of unity can be given by a
straightforward modification.

Let us show
\(
  f\ast g\in \widehat{\mathscr K}_t
\)
for $f$, $g\in \widehat{\mathscr K}_t$.
By induction and \eqref{eq:example} we may assume
that $f$ is of form
\begin{equation*}
   m'\left(1 + V_{i,b\ve}\right),
\end{equation*}
where $m'$ is a monomial with $u_{i,b}(m') = 1$, $u_{i,c}(m') = 0$ for
$c\neq b$, and that $g = E_i(m)$ is as \eqref{eq:form}.
By a direct calculation, we get
\begin{multline*}
   t^{-2d(m',m)} f\ast g - E_i(mm')
\\
=  \left(t^{2n} - 1 \right) m m'
     \prod_{a\neq b\ve^{-2}}
     \sum_{r_a=0}^{u_{i,a}(m)}
     t^{r_a(u_{i,a}(m)-r_a)}
     \begin{bmatrix}
       u_{i,a}(m) \\ r_a
     \end{bmatrix}_t
     \,
     V_{i,a\ve}^{r_a}
     \sum_{s=0}^{n-1} t^{s(n-s)}
     \begin{bmatrix}
       n-1 \\ s
     \end{bmatrix}_t
     V_{i,b\ve^{-1}}^{s+1}
\end{multline*}
where $n=u_{i,b\ve^{-2}}(m)$. If $n=0$, then the right hand side is
zero, so the assertion is obvious.
If $n\neq 0$, then we have
\begin{equation*}
   u_{i,a}\left(mm'V_{i,b\ve^{-1}}\right)
   = 
   \begin{cases}
     u_{i,b\ve^{-2}}(m) - 1 & \text{if $a = b\ve^{-2}$},
   \\
     u_{i,a}(m) & \text{otherwise}.
   \end{cases}
\end{equation*}
Therefore the above expression is equal to 
\(
   \left(t^{2n} - 1 \right)E_i\left(mm'V_{i,b\ve^{-1}}\right)
\).

Next we show the closedness of the image under the involution. By
\eqref{eq:antihom} and the above assertion, we may assume 
\(
   f = m'\left(1 + V_{i,b\ve}\right)
\)
as above. We further assume $m'$ does not contain $t$, $t^{-1}$. Then
we get
\begin{equation*}
   \overline{f} = t^{2d(m',m')} f.
\end{equation*}
This is contained in $\widehat{\mathscr K}_t$.

Now we can define $\tilde\ast$ and
$\setbox5=\hbox{A}\overline{\rule{0mm}{\ht5}\hspace*{\wd5}}$
on $\bfR_t$ so that
\begin{gather*}
   \chi_{\ve,t}(\overline{V}) 
     = \widehat\Pi\left(\overline{\widehat\chi_{\ve,t}(V)}\right)
     = \overline{\chi_{\ve,t}(V)},
\\
   \chi_{\ve,t}(V_1\tilde\ast V_2)
     = \chi_{\ve,t}(V_1)\tilde\ast\chi_{\ve,t}(V_2),
\end{gather*}
where we have assumed that $\ve$ is not a root of unity for the second 
equality.
By the above discussion together with \eqref{eq:use}, the right hand
sides are contained in $\mathscr K_t$, and therefore in the image of
$\chi_{\ve,t}$ by \thmref{thm:cons}(1).
Since $\chi_{\ve,t}$ is injective by \thmref{thm:cons}(1),
$\overline{V}$, $V_1\ast V_2$ are well-defined.
\end{proof}

\begin{Remark}
In this article, the existence of $\widehat\chi_{\ve,t}$ satisfying the
axioms is provided by a geometric theory of quiver varieties. But the
author conjectures that there exists purely combinatorial proof of the 
existence, independent of quiver varieties or the representation
theory of quantum loop algebras.
When $\mathfrak g$ is of type $A$ or $D$, such a combinatorial
construction is possible \cite{Na-AD}.
When $\mathfrak g$ is $E_6$, $E_7$, an explict construction of
$\widehat\chi_{\ve,t}$ is possible with the use of a computer.
\end{Remark}

\section{Graded and cyclic quiver varieties}\label{sec:quiver}

Suppose that a finite graph $(I,E)$ of type $ADE$ is given. The set
$I$ is the set of vertices, while $E$ is the set of edges.

Let $H$ be the set of pairs consisting of an edge together with its
orientation. For $h\in H$, we denote by $\vin(h)$ (resp.\ $\vout(h)$)
the incoming (resp.\ outgoing) vertex of $h$.
For $h\in H$ we denote by $\overline h$ the same edge as $h$ with the
reverse orientation.
We choose and fix a function $\varepsilon\colon H \to \C^*$ such that
$\varepsilon(h) + \varepsilon(\overline{h}) = 0$ for all $h\in H$.

Let $V$, $W$ be $I\times \C^*$-graded vector spaces such that
its $(i\times a)$-component, denoted by $V_i(a)$, is finite dimensional
and $0$ at most finitely many $i\times a$.
In what follows we consider only $I\times\C^*$-graded vector spaces
with this condition.
For an integer $n$, we define vector spaces by
\begin{equation}
\label{eq:LE}
\begin{gathered}[m]
  \HomL(V, W)^{[n]} \defeq
  \bigoplus_{i\in I, a\in\C^*}
    \Hom\left(V_i(a), W_i(a\ve^{n})\right),
\\
  \HomE(V, W)^{[n]} \defeq
  \bigoplus_{h\in H, a\in\C^*}
    \Hom\left(V_{\vout(h)}(a), W_{\vin(h)}(a\ve^{n})\right).
\end{gathered}
\end{equation}

If $V$ and $W$ are $I\times\C^*$-graded vector spaces as above, we
consider the vector spaces
\begin{equation}\label{def:bM}
  \bM \equiv \bM(V, W) \defeq
  \HomE(V, V)^{[-1]} \oplus \HomL(W, V)^{[-1]}
  \oplus \HomL(V, W)^{[-1]},
\end{equation}
where we use the notation $\bM$ unless we want to specify $V$, $W$.
The above three components for an element of $\bM$ is denoted by $B$,
$\alpha$, $\beta$ respectively.
({\bf NB}: In \cite{Na-qaff} $\alpha$ and $\beta$ were denoted by $i$, 
$j$ respectively.)
The
$\Hom\left(V_{\vout(h)}(a),V_{\vin(h)}(a\ve^{-1})\right)$-component of 
$B$ is denoted by $B_{h,a}$. Similarly, we denote by $\alpha_{i,a}$,
$\beta_{i,a}$ the components of $\alpha$, $\beta$.

We define a map $\mu\colon\bM\to \HomL(V,V)^{[-2]}$ by
\begin{equation*}
   \mu_{i,a}(B,\alpha,\beta)
   =  \sum_{\vin(h)=i} \ve(h)
      B_{h,a\ve^{-1}} B_{\overline{h},a} +
   \alpha_{i,a\ve^{-1}}\beta_{i,a},
\end{equation*}
where $\mu_{i,a}$ is the $(i,a)$-component of $\mu$.

Let $G_V \defeq \prod_{i,a} \GL(V_i(a))$. It acts on $\bM$ by
\begin{equation*}
  (B, \alpha, \beta) \mapsto g\cdot (B, \alpha, \beta)
  \defeq \left(g_{\vin(h),a\ve^{-1}} B_{h,a} g_{\vout(h),a}^{-1},\,
  g_{i,a\ve^{-1}}\alpha_{i,a},\,
  \beta_{i,a} g_{i,a}^{-1}\right).
\end{equation*}
The action preserves the subvariety $\mu^{-1}(0)$ in $\bM$.

\begin{Definition}\label{def:stable}
A point $(B, \alpha, \beta) \in \mu^{-1}(0)$ is said to be {\it stable\/} if 
the following condition holds:
\begin{itemize}
\item[] if an $I\times\C^*$-graded subspace $S$ of $V$ is
$B$-invariant and contained in $\Ker \beta$, then $S = 0$.
\end{itemize}
Let us denote by $\mu^{-1}(0)^{\operatorname{s}}$ the set of stable points.
\end{Definition}
Clearly, the stability condition is invariant under the action of
$G_V$. Hence we may say an orbit is stable or not.

We consider two kinds of quotient spaces of $\mu^{-1}(0)$:
\begin{equation*}
   \N_0(V,W) \defeq \mu^{-1}(0)\dslash G_V, \qquad
   \N(V,W) \defeq \mu^{-1}(0)^{\operatorname{s}}/G_V.
\end{equation*}
Here $\dslash$ is the affine algebro-geometric quotient, i.e., the
coordinate ring of $\N_0(V,W)$ is the ring of $G_V$-invariant
functions on $\mu^{-1}(0)$. In particular, it is an affine variety. It
is the set of closed $G_V$-orbits.
The second one is the set-theoretical quotient, but coincides with a
quotient in the geometric invariant theory (see \cite[\S3]{Na:1998}).
The action of $G_V$ on $\mu^{-1}(0)^{\operatorname{s}}$ is free thanks 
to the stability condition (\cite[3.10]{Na:1998}).
By a general theory, there exists a natural projective morphism
\begin{equation*}
   \pi\colon \N(V,W) \to \N_0(V,W).
\end{equation*}
(See \cite[3.18]{Na:1998}.)
The inverse image of $0$ under $\pi$ is denoted by $\NLa(V,W)$.
We call these varieties {\it cyclic quiver varieties\/} or {\it graded
quiver varieties}, according as $\ve$ is a root of unity or {\it not}.

Let $\Nreg(V,W)\subset\N_0(V,W)$ be a possibly empty open subset of
$\N_0(V,W)$ consisting of free $G_V$-orbits. It is known that
$\pi$ is isomorphism on $\pi^{-1}(\Nreg(V,W))$ \cite[3.24]{Na:1998}.
In particular, $\Nreg(V,W)$ is nonsingular and is pure dimensional.

A $G_V$-orbit though $(B,\alpha,\beta)$, considered as a point of
$\N(V,W)$ is denoted by $[B,\alpha,\beta]$.

We associate polynomials $e^{W}$, $e^{V}\in\widehat{\mathscr
Y}_t$ to graded vector spaces $V$, $W$ by
\begin{equation}\label{eq:rule}
   e^{W} = \prod_{i\in I,a\in\C^*} W_{i,a}^{\dim W_i(a)},
\quad
   e^{V} = \prod_{i\in I,a\in\C^*} V_{i,a}^{\dim V_i(a)}.
\end{equation}

Suppose that we have two $I\times\C^*$-graded vector spaces $V$, $V'$
such that $V_i(a) \subset V'_i(a)$ for all $i$, $a$. Then $\N_0(V,W)$
can be identified with a closed subvariety of $\N_0(V',W)$ by the
extension by $0$ to the complementary subspace (see
\cite[2.5.3]{Na-qaff}). We consider the limit
\begin{equation*}
   \N_0(\infty,W) \defeq \bigcup_{V} \N_0(V,W).
\end{equation*}
It is known that the above stabilizes at some $V$ (see \cite[2.6.3,
2.9.4]{Na-qaff}).
The complement $\N_0(V,W)\setminus\Nreg(V,W)$ consists of a finite
union of $\Nreg(V',W)$ for smaller $V'$'s \cite[3.27,
3.28]{Na:1998}. Therefore we have a decomposition
\begin{equation}\label{eq:stratum}
   \N_0(\infty,W) = \bigsqcup_{[V]} \N_0(V,W),
\end{equation}
where $[V]$ denotes the isomorphism class of $V$. The transversal
slice to each stratum was constructed in \cite[\S3.3]{Na-qaff}.
Using it, we can check
\begin{align}
&\begin{minipage}[m]{0.75\textwidth}
\noindent
If $\Nreg(V,W)\neq\emptyset$, then
$e^{V} e^{W}$ is {\it l\/}--dominant.
\end{minipage}\label{eq:l-dom}
\\
&\begin{minipage}[m]{0.75\textwidth}
\noindent
If $\Nreg(V,W)\subset\overline{\Nreg(V',W)}$, then
$e^{V'} \le e^V$.
\end{minipage}\label{eq:closure}
\end{align}

On the other hand, we consider the disjoint union for $\N(V,W)$:
\[
   \N(W) \defeq \bigsqcup_{[V]} \N(V,W).
\]
Note that there are no obvious morphisms between $\N(V,W)$ and
$\N(V',W)$ since the stability condition is not preserved under the
extension.
We have a morphism $\N(W)\to \N_0(\infty,W)$, still denoted by $\pi$.

The original quiver varieties \cite{Na:1994,Na:1998} are the special
case when $\ve = 1$ and $V_i(a) = W_i(a) = 0$ except $a=1$.
On the other hand, the above varieties $\N(W)$, $\N_0(\infty,W)$ are
fixed point set of the original quiver varieties with respect to a
semisimple element in a product of general linear groups. (See
\cite[\S4]{Na-qaff}.) In particular, it follows that $\N(V,W)$ is
nonsingular, since the corresponding original quiver variety is so.
This can be also checked directly.

Since the action is free, $V$ and $W$ can be considered as
$I\times\C^*$-graded vector bundles over $\N(V,W)$. We denote them by the
same notation.  We consider $\HomE(V,V)$, $\HomL(W,V)$, $\HomL(V, W)$
as vector bundles defined by the same formula as in \eqref{eq:LE}. By
the definition, $B$, $\alpha$, $\beta$ can be considered as sections of those
bundles.

We define a three term sequence of vector bundles over $\N(V,W)$ by
\begin{equation}
\label{eq:taut_cpx_fixed}
  C_{i,a}^\bullet(V,W):
  V_i(a\ve)
\xrightarrow{\sigma_{i,a}}
  \displaystyle{\bigoplus_{h:\vin(h)=i}}
     V_{\vout(h)}(a)
    \oplus W_i(a)
\xrightarrow{\tau_{i,a}}
  V_i(a\ve^{-1}),
\end{equation}
where
\begin{equation*}
  \sigma_{i,a} = \bigoplus_{\vin(h)=i} B_{\overline h,a\ve}
      \oplus \beta_{i,a\ve},
  \qquad
  \tau_{i,a} = \sum_{\vin(h)=i} \varepsilon(h) B_{h,a} + \alpha_{i,a}.
\end{equation*}
This is a complex thanks to the equation $\mu(B,\alpha,\beta) = 0$.
We assign the degree $0$ to the middle term.
By the stability condition, $\sigma_{i,a}$ is injective.

We define the rank of complex $C^\bullet$ by $\sum_p (-1)^p \rank
C^p$. Then we have
\begin{equation*}
   \rank C_{i,a}^\bullet(V,W) = u_{i,a}(e^{V}e^{W}).
\end{equation*}
We denote the right hand side by $u_{i,a}(V,W)$ for brevity.

We define a three term complex of vector bundles over
$\N(V^1,W^1)\times \N(V^2,W^2)$ by
\begin{equation}\label{eq:hecke_complex}
   \HomL(V^1,V^2)^{[0]}
   \xrightarrow{\sigma^{21}}
   \begin{matrix}
     \HomE(V^1, V^2)^{[-1]} \\
     \oplus \\
     \HomL(W^1, V^2)^{[-1]} \\
      \oplus \\
     \HomL(V^1, W^2)^{[-1]}
   \end{matrix}
   \xrightarrow{\tau^{21}}
   \HomL(V^1,V^2)^{[-2]},
\end{equation}
where
\begin{equation*}
\begin{split}
        \sigma^{21}(\xi) & = (B^2 \xi - \xi B^1) \oplus
         (-\xi \alpha^1) \oplus \beta^2 \xi, \\
        \tau^{21}(C\oplus I\oplus J)
         &= \varepsilon B^2 C + \varepsilon C B^1 + \alpha^2 J + I \beta^1.
\end{split}
\end{equation*}
We assign the degree $0$ to the middle term.
By the same argument as in \cite[3.10]{Na:1998}, $\sigma^{21}$ is
injective and $\tau^{21}$ is surjective. Thus the quotient
$\Ker\tau^{21}/\Ima\sigma^{21}$ is a vector bundle over
$\N(V^1,W^1)\times \N(V^2, W^2)$. Its rank is given by
\begin{equation}\label{eq:rank}
  d(e^{V^1}e^{W^1}, e^{V^2}e^{W^2}).
\end{equation}

If $V^1=V^2$, $W^1=W^2$, then the restriction of
$\Ker\tau^{21}/\Ima\sigma^{21}$ to the diagonal is isomorphic to the
tangent bundle of $\N(V,W)$ (see \cite[Proof of 4.1.4]{Na-qaff}). In
particular, we have
\begin{equation}\label{eq:dim}
    \dim\N(V,W) = d(e^{V}e^{W},e^{V}e^{W}).
\end{equation}

Let us give the definition of $\widehat\chi_{\ve,t}$.
We define $\widehat\chi_{\ve,t}$ for all standard modules
$M(P)$.
Since $\{ M(P)\}_P$ is a basis of $\operatorname{Rep}\Ule$, we
can extend it linearly to any finite dimensional $\Ule$-modules.

The relation between standard modules and graded/cyclic quiver
varieties is as follows (see \cite[\S13]{Na-qaff}):
Choose $W$ so that $e^W = e^P$, i.e., 
\[
  P_i(u) = \prod_a (1-au)^{\dim W_{i}(a)}.
\]
Then a standard module $M(P)$ is defined as $H_*(\NLa(W),\C)$, which
is equipped with a structure of $\Ule$-module by the convolution
product. Moreover, its {\it l\/}--weight space $M(P)_Q$ is
\[
  \bigoplus_{V: e^Ve^W = e^Q} H_*(\NLa(V,W),\C).
\]
Here $H_k(\ ,\C)$ denotes the Borel-Moore homology with complex
coefficients.
If $\ve$ is not a root of unity, then $V$ is determined from $Q$. So
the above has only one summand.

Let
\begin{equation}\label{eq:geomdef}
   \widehat\chi_{\ve,t}(M(P))
   \defeq \sum_{[V]} (-t)^{k} \dim H_k(\NLa(V,W),\C)\, e^{V}e^{W}.
\end{equation}
Since $H_k(\NLa(V,W),\C)$ vanishes for odd $k$ \cite[\S7]{Na-qaff}, we
may replace $(-t)^k$ by $t^k$. In particular, it is clear that
coefficients of $\widehat\chi_{\ve,t}(M(P))$ are polynomials in $t^2$
with positive coefficients.

In subsequent sections we prove that the above $\widehat\chi_{\ve,t}$
satisfies the axioms.
By definition, it is clear that $\widehat\chi_{\ve,t}$ satisfies Axiom~1.

Remark that \corref{cor:coincide} directly follows from this geometric
definition (\cite[13.4.5]{Na-qaff}).

We give a simple consequence of the definition:
\begin{Proposition}\label{prop:nodom}
Assume $\ve$ is not a root of unity.
Suppose that all roots of $P_i(u) = 0$ have the same value
\textup(e.g., $P_i(u) = 1 - au$, $P_j(u) = 1$ for $j\neq i$
for some $i$\textup).
Then $M(P)$ has no {\it l\/}--dominant term other than $e^P$.
\end{Proposition}

This was proved in \cite[Corollary 4.5]{FM}. But we give a geometric
proof.

\begin{proof}
Take $W$ so that $e^W = e^P$. It is enough to show that $u_{i,a}(V,W)
< 0$ for some $i$, $a$ if $\N(V,W)\neq\emptyset$ and $V\neq 0$.

By the assumption, there is a nonzero constant $a$ such that $W_i(b) = 
0$ for all $i$, $b\neq a$. By the stability condition, we have
\(
   V_i(b) = 0
\)
if $b\neq a\ve^n$ for some $n\in\Z_{>0}$. Let $n_0$ be the maximum of
such $n$, and suppose $V_i(a\ve^{n_0})\neq 0$. Since
$W_i(a\ve^{n_0+1}) = V_i(a\ve^{n_0+1}) = V_i(a\ve^{n_0+2}) = 0$, we
have
\begin{equation*}
   u_{i,a\ve^{n_0+1}}(V,W) = \rank C_{i,a\ve^{n_0+1}}^\bullet(V,W) < 0.
\end{equation*}
\end{proof}

\section{Proof of Axiom~2: Analog of the Weyl group invariance}
\label{sec:proof1}


For a complex algebraic variety $X$, let $e(X;x,y)$ denote the virtual
Hodge polynomial defined by Danilov-Khovanskii~\cite{DK} using a mixed
Hodge strucuture of Deligne~\cite{De}. It has the following properties.
\begin{enumerate}
\item $e(X;x,y)$ is a polynomial in $x$, $y$ with integer coefficients.
\item If $X$ is a nonsingular projective variety, then $e(X;x,y) =
  \sum_{p,q} (-1)^{p+q} h^{p,q}(X)x^p y^q$, where the $h^{p,q}(X)$ are the
  Hodge numbers of $X$.
\item If $Y$ is a closed subvariety in $X$, then $e(X;x,y) = e(Y;x,y) +
e(X\setminus Y;x,y)$.
\item If $f\colon Y\to X$ is a fiber bundle with fiber $F$ which is
  locally trivial in the Zarisky topology, then $e(Y;x,y) =
  e(X;x,y)e(F;x,y)$.
\end{enumerate}

We define the virtual Poincar\'e polynomial of $X$ by $p_t(X) \defeq
e(X;t,t)$. (In fact, this reduction does not loose any
information. The following argument shows that $e(X;x,y)$ appearing
here is a polynomial in $xy$.)
The actual Poincar\'e polynomial is defined as
\begin{equation*}
  P_t(X) = \sum_{k=0}^{2\dim X} (-t)^k \dim H_k(X,\C),
\end{equation*}
where $H_k(X,\C)$ is the Borel-Moore homology of $X$ with complex
coefficients.

\begin{Remark}
In stead of virtual Poincar\'e polynomials, we can use numbers of
rational points in the following argument, if we define graded/cyclic
varieties over an algebraic closure of a finite field $\mathbf k$. As
a consequence, those numbers are special values of ``computable''
polynomials $P(t)$ at $t=\sqrt{\#\mathbf k}$.
\end{Remark}

\begin{Lemma}\label{lem:vir}
The virtual Poincar\'e polynomial of $\NLa(V,W)$ is equal to the actual
Poincar\'e polynomial. Moreover, it is a polynomial in $t^2$.
The same holds for $\N(V,W)$.
\end{Lemma}

\begin{proof}
In \cite[\S7]{Na-qaff} we have shown that $\NLa(V,W)$ has
an partition into locally closed subvarieties $X_1,\dots, X_n$ with
the following properties:
\begin{enumerate}
\item $X_1\cup X_2\cup\dots \cup X_i$ is closed in $\NLa(V,W)$ for
each $i$.
\item each $X_i$ is a vector bundle over a nonsingular projective
variety whose homology groups vanish in odd degrees.
\end{enumerate}
A partition satisfying the property (1) is called an {\it
$\alpha$-partition}.
(More precisely, it was shown in [loc.\ cit.] that $X_i$ is a fiber
bundle with an affine space fiber over the base with the above
property. The above statement was shown in \cite{Na-tensor}.)
  
By the long exact sequence in homology groups, we have
$P_t(\NLa(V,W)) = \sum_i P_t(X_i)$. On the other hand, by the
property of the virtual Poincar\'e polynomial, we have
$p_t(\NLa(V,W)) = \sum_i p_t(X_i)$. Since $X_i$ satisfies the
required properties in the statement, it follows that $\NLa(V,W)$
satisfies the same property.
  
We have an $\alpha$-partition with the property (2) also for
$\N(V,W)$, so we have the same assertion.
\end{proof}

Recall the complex \eqref{eq:taut_cpx_fixed}. For a $\C^*$-tuple of
nonnegative integers $(n_a)\in\Z_{\ge 0}^{\C^*}$, let
\begin{equation*}
  \N_{i;(n_a)}(V,W)
  \defeq 
  \left\{ [B,\alpha,\beta]\in \N(V,W) \Biggm| \text{$\codim_{V_i(\ve^{-1}a)}
        \Ima \tau_{i,a}
        = n_a$ for each $a\in\C^*$}\right\}.
\end{equation*}
This is a locally closed subset of $\N(V,W)$.
We also set
\begin{equation*}
  \NLa_{i;(n_a)}(V,W) \defeq \N_{i;(n_a)}(V,W)\cap\NLa(V,W).
\end{equation*}
We have partitions
\begin{equation*}
   \N(V,W) = \bigsqcup_{(n_a)} \N_{i;(n_a)}(V,W), \qquad
   \NLa(V,W) = \bigsqcup_{(n_a)} \NLa_{i;(n_a)}(V,W).
\end{equation*}

Let $Q_{i,a}(V,W)$ be the middle cohomology of the
complex $C_{i,a}^\bullet(V,W)$ \eqref{eq:taut_cpx_fixed}, i.e.,
\begin{equation*}
   Q_{i,a}(V,W) \defeq \Ker\tau_{i,a}/\Ima\sigma_{i,a}.
\end{equation*}
Over each stratum $\N_{i;(n_a)}(V,W)$ it defines a vector bundle. In
particular, over the open stratum $\N_{i;(0)}(V,W)$, i.e., points
where $\tau_{i,a}$ is surjective for all $i$, its rank is equal to
\begin{equation}\label{eq:rankC}
  \rank C^\bullet_{i,a}(V,W) = u_{i,a}(V,W).
\end{equation}

Suppose that a point $[B,\alpha,\beta]\in\N_{i;(n_a)}(V,W)$ is given.
We define a new graded vector space $V'$ by $V'_i(\ve^{-1}a) \defeq
\Ima\tau_{i,a}$. The restriction of $(B,i,j)$ to $V'$ also satisfies
the equation $\mu = 0$ and the stability condition. Therefore it
defines a point in $\N(V',W)$.  It is clear that this construction
defines a map
\begin{equation}
\label{eq:projection}
   p\colon \N_{i;(n_a)}(V,W)\to \N_{i;(0)}(V',W).
\end{equation}

Let $G(n_a,Q_{i,\ve^{-2}a}(V',W)|_{\N_{i;(0)}(V',W)})$ denote the
Grassmann bundle of $n_a$-planes in the vector bundle obtained by
restricting $Q_{i,\ve^{-2}a}(V',W)$ to $\N_{i;(0)}(V',W)$. Let
\[
  {\displaystyle \prod_a}
     G(n_a,Q_{i,\ve^{-2}a}(V',W)|_{\N_{i;(0)}(V',W)})
\]
be their fiber product over $\N_{i;(0)}(V',W)$.
By \cite[5.5.2]{Na-qaff} there exists a commutative diagram
\begin{equation*}
\begin{CD}
   {\displaystyle \prod_a}
     G(n_a,Q_{i,\ve^{-2}a}(V',W)|_{\N_{i;(0)}(V',W)})
      @>\pi>> \N_{k;(0)}(V',W) \\
   @VV{\cong}V @| \\
   \N_{i;(n_a)}(V,W) @>p>> \N_{i;(0)}(V',W),
\end{CD}
\end{equation*}
where $\pi$ is the natural projection. (The assumption $\ve\neq\pm 1$
there was unnecessary. See \secref{sec:pm 1}.)

Since the projection~\eqref{eq:projection} factors through $\pi$, it
induces
\[
   p\colon \NLa_{i;(n_a)}(V,W)\to \NLa_{i;(0)}(V',W).
\]
Therefore we have
\begin{equation*}
\begin{split}
   & 
   p_t\left(\NLa_{i;(n_a)}(V,W)\right)
\\
   =\; &\prod_{a}
   t^{n_a(\rank C^\bullet_{i,\ve^{-2}a}(V',W)-n_a)}
   \begin{bmatrix}
     \rank C^\bullet_{i,\ve^{-2}a}(V',W)
     \\
     n_a
   \end{bmatrix}_t
   p_t\left(\NLa_{i;(0)}(V',W)\right).
\end{split}
\end{equation*}
Using \eqref{eq:rankC}, we get
\begin{equation*}
\begin{split}
  \widehat\chi_{\ve,t}(M(P)) = \sum_{[V']}\;
   & p_t\left(\NLa_{i;(0)}(V',W)\right) e^{V'}e^{W}
\\
   & \times \prod_a
     t^{r_a(u_{i,a}(V',W)-r_a)}
     \sum_{r_a=0}^{u_{i,a}(V',W)}
     \begin{bmatrix}
       u_{i,a}(V',W) \\ r_a
     \end{bmatrix}_t V_{i,a\ve}^{r_a}.
\end{split}
\end{equation*}
This shows that $\widehat\chi_{\ve,t}(M(P))$ is contained
$\widehat{\mathscr K}_t$. This completes the proof of
Axiom~2.

As we promised, we give an algorithm computing
$\widehat\chi_{\ve,t}(M)$ for $M$ as in \thmref{thm:cons}(2).
Although we will explain it only when $M$ is a standard module $M(P)$,
a modification to the general case is straightforward, if we interpret
$p_t(\NLa_{i;(n_a)}(V,W))$ suitably.
Moreover, we use the Grassmann bundle \eqref{eq:projection} instead of 
the condition $\widehat\chi_{\ve,t}(M)\in\widehat{\mathscr K}_t$.

We ``compute'' the virtual Poincar\'e polynomials
$p_t(\NLa_{i;(n_a)}(V,W))$ by induction.
The first step of the induction is the case $V = 0$.
In this case, $\NLa(0,W)$ is a single point, so $p_t(\NLa(0,W)) = 1$.
Moreover, $\NLa(0,W) = \NLa_{i;(0)}(0,W)$ for all $i$.

Suppose that we already ``compute'' all
$p_t(\NLa_{i;(n'_\lambda)}(V',W))$ for $\dim V' < \dim V$.
If $(n_a)\neq (0)$, then by the Grassmann bundle
\eqref{eq:projection} and the induction hypothesis, we can ``compute'' 
$p_t(\NLa_{i;(n_a)}(V,W))$.
By the assumption in \thmref{thm:cons}(2), $e^{V}e^{W}$ is {\it not\/}
{\it l\/}--dominant, hence there exists $i$, $a$ such that
$u_{i,a}(V,W) < 0$.
Then $\NLa_{i;(0)}(V,W)$ is the empty set by
\cite[5.5.5]{Na-qaff}. Therefore we have
\begin{equation*}
   p_t\left(\NLa(V,W)\right)
   = \sum_{(n_a)\neq (0)} p_t\left(\NLa_{i;(n_a)}(V,W)\right).
\end{equation*}
We have already ``computed'' the right hand side. We, of course, have
\begin{equation*}
   p_t\left(\NLa_{i;(0)}(V,W)\right) = 0.
\end{equation*}
For $j\neq i$, we have
\begin{equation*}
   p_t\left(\NLa_{j;(0)}(V,W)\right)
   = p_t\left(\NLa(V,W)\right)
   - \sum_{(n_a)\neq (0)} p_t\left(\NLa_{j;(n_a)}(V,W)\right).
\end{equation*}
The right hand side is already ``computed''.

\begin{Remark}
(1) Note that the above argument shows that $p_t(\NLa(V,W))$ is a
polynomial in $t^2$ without appealing to \cite[\S7]{Na-qaff} as in
\lemref{lem:vir}.
In fact, nonsingular quasi-projective varieties appearing in
\lemref{lem:vir} are examples of graded quiver varieties such that the
above argument can be applied, i.e., the corresponding standard
modules satisfy the condition in \thmref{thm:cons}(2). Therefore, the
above gives a new proof of the vanishing of odd homology groups.
  
(2) If the reader carefully compares our algorithm with
Frenkel-Mukhin's one \cite{FM}, he/she finds a difference. The {\it
coloring\/} $s_i$ of a monomial $m = e^{V}e^{W}$ is
\[
   \sum_{(n_a)\neq (0)}
   p_{t=1}\left(\NLa_{i;(n_a)}(V,W)\right)
\]
in our algorithm. This might be possibly {\it negative\/} integer,
while it is assumed to be nonnegative in [loc.\ cit.]. Therefore, we
must modify their definition of the {\it admissibility\/} of a
monomial $m$.
Let us consider all values $s_i$ such that $m$ is not
$i$--dominant. We say $m$ is {\it admissible\/} if all values are the
same.
In our case, $s_i$ is
\(
   p_{t=1}\left(\NLa(V,W)\right)
\)
if $m$ is not $i$--dominant, hence $\NLa_{i;(0)}(V,W)=\emptyset$.
Therefore it is independent of $i$.
\end{Remark}

\section{Proof of Axiom~3: Multiplicative property}
\label{sec:proof2}

Since \cite{Na-hom} it has been known that Betti numbers of arbitrary
quiver varieties are determined by special cases corresponding to
fundamental weights. We will use the same idea in this section.

Let $W^1$, $W^2$, $W$ be $I\times\C^*$-graded vector spaces such that 
$W_i(a) = W^1_i(a)\oplus W^2_i(a)$ for $i\in I$, $a\in\C^*$.
Let $P^1$, $P^2$ be $I$-tuples of polynomials corresponding to $W^1$,
$W^2$. Then $W$ corresponds to $P^1P^2$.

We define a map $\N(W^1)\times \N(W^2)\to \N(W)$ by
\begin{equation}\label{eq:directsum}
   \left([B^1,\alpha^1,\beta^1], [B^2,\alpha^2,\beta^2]\right)
   \longmapsto
   [B^1\oplus B^2, \alpha^1\oplus\alpha^2, \beta^1\oplus \beta^2].
\end{equation}

We define a $\prod_{i,a} \GL(W_i(a))$-action on $\N(W)$ by
\begin{equation*}
   s \ast [B,\alpha,\beta] \defeq [B, \alpha s^{-1}, s\beta].
\end{equation*}

We define a one-parameter subgroup $\lambda\colon\C^*\to \prod_{i,a}
\GL(W_i(a))$ by
\begin{equation*}
     \lambda(t) = \bigoplus_{i,a}
        \id_{W^1_i(a)} \oplus\, t\id_{W^2_i(a)}
.
\end{equation*}

Then \eqref{eq:directsum} is a closed embedding and the fixed point
set $\N(W)^{\lambda(\C^*)}$ is its image by \cite[3.2]{Na-tensor}. We
identify $\N(W^1)\times \N(W^2)$ with its image hereafter.
The fixed point set $\N(V,W)^{\lambda(\C^*)}$ is the union of
$\N(V^1,W^1)\times\N(V^2,W^2)$ with $V\cong V^1\oplus V^2$.

Let
\begin{multline*}
   \Zm(V^1,W^1;V^2,W^2) \defeq
\\
   \left\{ [B,\alpha,\beta]\in\N(W) \left|\,
       \lim_{t\to 0}\lambda(t)\ast[B,\alpha,\beta]\in
       \N(V^1,W^1)\times\N(V^2,W^2)
     \right\}\right..
\end{multline*}
We also define $\Zl(V^1,W^1;V^2,W^2)$ by replacing
$\N(V^1,W^1)\times\N(V^2,W^2)$ by $\NLa(V^1,W^1)\times\NLa(V^2,W^2)$.
These are studied in \cite{Na-tensor}. These are nonsingular locally
closed subvarieties of $\N(W)$ [loc.\ cit., 3.7].

Let $\Zm(W^1;W^2)$, $\Zl(W^1,W^2)$ be their union over $[V^1]$,
$[V^2]$ respectively. These are closed subvarieties of
$\N(W)$ [loc.\ cit., 3.6].
By [loc.\ cit., 3.7, 3.13], the partition
\[
   \Zm(W^1;W^2) = \bigsqcup_{[V^1],[V^2]} \Zm(V^1,W^1;V^2,W^2) 
\]
is an $\alpha$-partition such that each stratum $\Zm(V^1,W^1;V^2,W^2)$
is isomorphic to the total space of the vector bundle
$\Ker\tau^{21}/\Ima\sigma^{21}$ over $\N(V^1,W^1)\times\N(V^2,W^2)$
in \eqref{eq:hecke_complex}. (More precisely, we restrict the result
of [loc.\ cit.] to the fixed point set.)

Similarly
\[
   \Zl(W^1;W^2) = \bigsqcup_{V^1,V^2} \Zl(V^1,W^1;V^2,W^2) 
\]
is an $\alpha$-partition such that each stratum $\Zm(V^1,W^1;V^2,W^2)$
is isomorphic to the restriction of $\Ker\tau^{21}/\Ima\sigma^{21}$ to 
$\NLa(V^1,W^1)\times\NLa(V^2,W^2)$.

\begin{Proposition}\label{prop:mul}
\textup{(1)} The virtual Poincar\'e polynomial of $\Zl(W^1;W^2)$
\textup(more precisely that of each connected component of
$\Zl(W^1;W^2)$\textup) is equal to its actual Poincar\'e
polynomial. Moreover, it is a polynomial in $t^2$.
The same holds for $\Zm(W^1;W^2)$.

\textup{(2)} We have
\begin{equation}\label{eq:generalize}
  \widehat\chi_{\ve,t}(M(P^1)) \ast \widehat\chi_{\ve,t}(M(P^2)) =
  \sum_{[V^1],[V^2]}
  P_t(\Zl(V^1,W^1;V^2,W^2))\; e^{V^1} e^{V^2} e^{W^1} e^{W^2}.
\end{equation}

\textup{(3)} The above expression is contained in $\widehat{\mathscr
K}_t$.
\end{Proposition}

\begin{proof}
(1) This can be shown exactly as in \lemref{lem:vir}.
  
(2) The rank of the vector bundle $\Ker\tau^{21}/\Ima\sigma^{21}$ is
equal to $d(e^{W^1}e^{V^1},e^{W^2}e^{V^2})$ (see \eqref{eq:rank}). By
the property of virtual Poincar\'e polynomials, we get the assertion.

(3) Exactly the same as \secref{sec:proof1}.
\end{proof}

Axiom~3 follows from above and the following assertion proved in
\cite[6.12]{Na-tensor}:
\begin{equation*}
    \Zl(W^1;W^2) = \NLa(W)
\end{equation*}
under the condition \eqref{eq:Z}.

As we promised, we give
\begin{proof}[A different proof of \thmref{thm:bar}]
We only prove the second statement. In fact, it is not difficult to
show that the following argument also implies the first statement.
  
We will prove that our $\widehat\chi_{\ve,t}$ satisfies Axiom~4 in the
next section. Therefore, it is enough to check the assertion for
$\widehat\chi_{\ve,t}$ given by the geometric
definition~\eqref{eq:geomdef}.

By \thmref{thm:cons}(1) and \propref{prop:mul}, we get the statement
for the multiplication.

Similarly, for the proof of the statement for the involution, it is
enough to show
\begin{equation*}
   \overline{\widehat\chi_{\ve,t}(M(P))} \in \widehat{\mathscr K}_t.
\end{equation*}
This follows from
\begin{equation}\label{eq:Claim}
   \overline{\widehat\chi_{\ve,t}(M(P))}
   =   \displaystyle\sum_{[V]} P_t(\N(V,W))\; e^{V}e^{W}.
\end{equation}

In fact, the same argument as in the proof of \secref{sec:proof1}
shows that the right hand side is contained in $\widehat{\mathscr
K}_t$.

Let us prove \eqref{eq:Claim}.
Since $\NLa(V,W)$ is homotopic to $\N(V,W)$ \cite[4.1.2]{Na-qaff}, its
usual homology group is isomorphic to that of $\N(V,W)$. Since
$\NLa(V,W)$ is compact, the usual homology group is isomorphic to the
Borel-Moore homology. Therefore, the Poincar\'e duality for $\N(V,W)$,
which is applicable since $\N(V,W)$ is nonsingular, implies
\begin{equation*}
   t^{2\dim \N(V,W)}
      P_{1/t}\left(\NLa(V,W)\right)
   =  P_{t}\left(\N(V,W)\right).
\end{equation*}
Since $\dim\N(V,W) = d(e^{V}e^{W},e^{V}e^{W})$
(see \eqref{eq:dim}), we get \eqref{eq:Claim}.
\end{proof}

\section{Proof of Axiom~4: Roots of unity}\label{sec:proof3}

In this section, we use a $\C^*$-action on $\N(V,W)$ to calculate
Betti numbers. This idea was originally appeared in \cite{Na-hom} and
\cite[\S5]{Na:1994}.

We assume that $\ve$ is a primitive $s$-th root of unity
($s\in\Z_{>0}$).
%
%

We may assume $\alpha = 1$ in the setting of Axiom~4.
We consider $V$, $W$ as $I\times(\Z/s\Z)$-graded vector
spaces.

We define a $\C^*$-action on $\bM(V,W)$ by
\[
   t\star (B,\alpha,\beta) = (tB,t\alpha s(t)^{-1},t s(t)\beta),
   \qquad (t\in \C^*),
\]
where $s(t)\in \prod \GL(W_i(a))$ is defined by
\begin{equation*}
   s(t) = \bigoplus_{i\in I,\; 0\le n < s} t^n\id_{W_i(\ve^n)}.
\end{equation*}
It preserves the equation $\mu(B,\alpha,\beta) = 0$ and commutes with
the action of $G_V$. Therefore it induces an action on the affine
cyclic quiver variety $\N_0(V,W)$. The action preserves the stability
condition. Therefore it induces action on $\N(V,W)$. These induced
actions are also denoted by $\star$. The map
$\pi\colon\N(V,W)\to\N_0(V,W)$ is equivariant.

\begin{Lemma}
Let $[B,\alpha,\beta]\in\N(V,W)$. Consider the flow $t\star
[B,\alpha,\beta]$ for $t\in\C^*$. It has a limit when $t\to 0$.
\end{Lemma}

\begin{proof}
By a generaly theory, it is enough to show that
$t\star[B,\alpha,\beta]$ stays in a compact set. Since $\pi$ is
proper, it is enough to show that $\pi(t\star[B,\alpha,\beta]) =
t\star\pi([B,\alpha,\beta])$ stays in a compact set.

By \cite{Lu:qv} the coordinate ring of $\N_0(V,W)$ is generated by
functions of forms
\begin{equation*}
  \langle  \chi, \beta_{i,\ve^{n+1}}
  B_{h_1,\ve^{n+2}}\dots B_{h_N,\ve^{n+N+1}}
  \alpha_{j,\ve^{n+N+2}}\rangle
\end{equation*}
where $\chi$ is a linear form on $\Hom\left(W_{i,\ve^n},
  W_{j,\ve^{n+N+2}}\right)$, and $i = \vin(h_1)$, $\vout(h_1) =
\vin(h_2)$, \dots, $\vout(h_N) = j$. By the $\C^*$-action, this
function is multiplied by
\begin{equation*}
   t^{n+1} t^N t^{-r+1} = t^{n+N+2-r}
\end{equation*}
where we assume $0\le n < s$ and $r$ is the integer such that
$0\le r < s$ and $r \equiv n+N+2\mod s$. Then $n+N+2-r$ is
nonnegative. Therefore $t\star\pi([B,\alpha,\beta])$ stays in a
compact set for any $[B,\alpha,\beta]$.
\end{proof}

We want to identify a fixed point set in $\N(V,W)$ with some quiver
variety $\N(V_q,W_q)$ defined for $q$ which is {\it not\/} a root of
unity. We first explain a morphism from $\N(V_q,W_q)$ to $\N(V,W)$.
Vector spaces corresponding $\N(V_q,W_q)$ are $I\times\Z$-graded
vector spaces.
Suppose that $V_q, W_q$ are $I\times\Z$-graded vector space such that
$(W_q)_i(q^k) = 0$ unless $0\le k < s$ (no condition for $V_q$).
We consider $W_q$ as an $I\times(\Z/s\Z)$-graded vector space simply
identifying $\Z/s\Z$ with $\{ 0, 1, \dots, s-1\}$. Let us denote by
$W$ the resulting $I\times(\Z/s\Z)$-graded vector space.
We define an $I\times(\Z/s\Z)$-graded vector space $V$ by
\[
   V_i(\ve^n) \defeq \bigoplus_{k\not\equiv n\mod s} (V_q)_i(q^k).
\]
If a point in $\bM(V_q,W_q)$ is given, it defines a point in
$\bM(V,W)$ in obvious way. The map $\bM(V_q,W_q)\to \bM(V,W)$
preserves the equation $\mu = 0$ and the stability condition. It is
equivariant under the $G_{V_q}$ action, where $G_{V_q}\to G_V$ is an
obvious homomorphism. Therefore, we have a morphism
\begin{equation}
\label{eq:fixed}
   \N(V_q,W_q) \to \N(V,W).
\end{equation}
Note that $W_q$ is uniquely determined by $W$, while $V_q$ is {\it
not\/} determined from $V$.

\begin{Lemma}
A point $[B,\alpha,\beta]\in\N(V,W)$ is fixed by the $\C^*$-action if
and only if it is contained in the image of
\eqref{eq:fixed} for some $V_q$.
Moreover, the map \eqref{eq:fixed} is a closed embedding.
\end{Lemma}

\begin{proof}
Fix a representatitve $(B,\alpha,\beta)$ of $[B,\alpha,\beta]$.
Then $[B,\alpha,\beta]$ is a fixed point if and only if there
exists $\lambda(t)\in G_V$ such that
\begin{equation*}
   t\star (B,\alpha,\beta) = \lambda(t)^{-1}\cdot (B,\alpha,\beta).
\end{equation*}
Such $\lambda(t)$ is unique since the action of $G_V$ is free. In
particular, $\lambda\colon\C^*\to G_V$ is a group homomorphism.

Let $V_i(\ve^n)[k]$ be the weight space of $V_i(\ve^n)$ with
eigenvalue $t^k$. The above equation means that
\begin{gather*}
   B_{h,\ve^{n+1}}\left(V_{\vout(h)}(\ve^{n+1})[{k+1}]\right)
   \subset V_{\vin(h)}(\ve^{n})[{k}],
\qquad
   \alpha_{i,\ve^{n+1}}(W_i(\ve^{n+1})) \subset V_i(\ve^{n})[{n}],
\\
   \beta_{i,\ve^n}\left(V_i(\ve^n)[k]\right) = 0\quad
   \text{if $k\neq n$}.
\end{gather*}
Let us define an $I\times\C^*$-graded subspace $S$ of $V$ by
\begin{equation*}
    S_i(\ve^n) \defeq \bigoplus_{k\not\equiv n\mod s} V_i(\ve^n)[k].
\end{equation*}
The above equations imply that $S$ is contained in $\Ker\beta$ and
$B$-invariant. Therefore $S = 0$ by the stability condition.
This means that $[B,\alpha,\beta]$ is in the image of \eqref{eq:fixed} 
if we set
\begin{equation*}
   (V_q)_i(q^k) \defeq V_i(\ve^k)[k].
\end{equation*}

Conversely, a point in the image is a fixed point. Since $\lambda$ is
unique, the map \eqref{eq:fixed} is injective.

Let us consider the differential of \eqref{eq:fixed}. The tangent
space of $\N(V,W)$ at $[B,i,j]$ is the middle cohomology groups of
the complex
\begin{equation}\label{eq:tangent}
   \HomL(V,V)^{[0]}
   \xrightarrow{\sigma^{21}}
   \begin{matrix}
     \HomE(V, V)^{[-1]} \\
     \oplus \\
     \HomL(W, V)^{[-1]} \\
      \oplus \\
     \HomL(V, W)^{[-1]}
   \end{matrix}
   \xrightarrow{\tau^{21}}
   \HomL(V,V)^{[-2]}.
\end{equation}
Similarly the tangent space of $\N(V_q,W_q)$ is the middle cohomology
of a complex with $V$, $W$ are replaced by $V_q$, $W_q$.
We have a natural morphism between the complexes so that 
the induced map between cohomology groups is the differential of
\eqref{eq:fixed}. It is not difficult to show the injectivity by using 
the stability condition.
\end{proof}

Let us consider the tangent space $T$ of $\N(V,W)$ at
$[B,\alpha,\beta]\in\N(V_q,W_q)\subset\N(V,W)$, which is the middle
cohomology of \eqref{eq:tangent}. Let $V = \bigoplus_k V[k]$ be the
weight space decomposition as in the proof of the above lemma. The
tangent space $T$ has a weight decomposition $T = \bigoplus_k T[k]$,
where $T[k]$ is the middle cohomology of
{\footnotesize
\begin{equation*}
   \bigoplus_n \HomL(V[n],V[{n+k}])^{[0]}
   \xrightarrow{\sigma^{21}}
   \begin{matrix}
     \bigoplus_n \HomE(V[n], V[{n+k-1}])^{[-1]} \\
     \oplus \\
     \bigoplus_n \HomL(W[n], V[{n+k-1}])^{[-1]} \\
      \oplus \\
     \bigoplus_n \HomL(V[n], W[{n+k-1}])^{[-1]}
   \end{matrix}
   \xrightarrow{\tau^{21}}
   \bigoplus_n \HomL(V[n],V[{n+k-2}])^{[-2]},
\end{equation*}}

\noindent
where $W[n] = W(\ve^n)$ if $0\le n < s$ and $0$ otherwise.
The rank of the complex is equal to
\[
\begin{cases}
  d_q(e^{V_q} e^{W_q}, e^{V_q} e^{W_q}[k])
    & \text{if $k\equiv 0\mod s$},
\\
  0 & \text{otherwise}.
\end{cases}
\]
Here $e^{V_q} e^{W_q}[k]$ is defined as in \eqref{eq:m_shift}.

We consider the Bialynicki-Birula decomposition of $\N(V,W)$:
\begin{equation*}
\begin{split}
   & \N(V,W) = \bigsqcup_{[V_q]} S(V_q,W_q),
\\
   & \qquad S(V_q,W_q) \defeq
   \left\{ x\in \N(V,W) \left|\; \lim_{t\to 0} t\star x
       \in \N(V_q,W_q) \right\}\right..
\end{split}
\end{equation*}
By a general theory, each $S(V_q,W_q)$ is a locally closed subvariety
of $\N(V,W)$, and the natural map $S(V_q,W_q)\to \N(V_q,W_q)$ is a
fiber bundle whose fiber is an affine space of dimension equal to
$\sum_{k > 0} \dim T[k]$. By the above formula, it is equal to
\[
   \sum_{k > 0}
     d_q(e^{V_q}e^{W_q}, e^{V_q}e^{W_q}[ks]).
\]
We write this number by $D^+(e^{V_q}e^{W_q})$.

By the property of virtual Poincar\'e polynomials, we have
\begin{equation*}
   P_t(\N(V,W)) = \sum_{[V_q]} t^{2D^+(e^{V_q}e^{W_q})}
                    P_t(\N(V_q,W_q)).
\end{equation*}
(Recall that the virtual Poincar\'e polynomials coincide with the
actual Poincar\'e polynomials for these varieties.)
Combining with an argument in the proof of \eqref{eq:Claim}, we have
\begin{equation*}
\begin{split}
   P_t(\NLa(V,W)) &
   = t^{2d(e^{V}e^{W},e^{V}e^{W})}
               P_{1/t}(\N(V,W))
\\
   & = \sum_{[V_q]} t^{2d(e^{V}e^{W},e^{V}e^{W})
           - 2D^+(e^{V_q}e^{W_q})} P_{1/t}(\N(V_q,W_q))
\\
   & = \sum_{[V_q]} t^{2d(e^{V}e^{W},e^{V}e^{W})
        - 2D^+(e^{V_q}e^{W_q})
        - 2d_q(e^{V_q}e^{W_q},e^{V_q}e^{W_q})}
      P_{t}(\NLa(V_q,W_q)).
\end{split}
\end{equation*}
Since
\begin{equation*}
   d(e^{V}e^{W},e^{V}e^{W})
   = \dim T = \sum_k \dim T[k]
   = \sum_{k} d_q(e^{V_q}e^{W_q},e^{V_q}e^{W_q}[ks]),
\end{equation*}
we have
\begin{equation*}
\begin{split}
   & d(e^{V}e^{W},e^{V}e^{W})
        - D^+(e^{V_q}e^{W_q})
        - d_q(e^{V_q}e^{W_q},e^{V_q}e^{W_q})
\\
   =\; & \sum_{k < 0}
        d_q(e^{V_q}e^{W_q},e^{V_q}e^{W_q}[ks])
   = D^-(e^{V_q}e^{W_q}).
\end{split}
\end{equation*}
Thus we have checked Axiom~4.

\begin{Remark}
When $\ve = 1$, there is a different $\C^*$-action so that the index
$D^-(m)$ can be read off from $a_m(t)$. See \cite[\S7]{Na-ann}.
\end{Remark}

\section{Perverse sheaves on graded/cyclic quiver varieties}
\label{sec:perverse}

The following is the main result of this article:
\begin{Theorem}\label{thm:character}
\textup{(1)} There exists a unique base $\{ L(P)\}$ of $\bfR_t$ such
that
\begin{equation*}
  \overline{L(P)} = L(P), \qquad
    L(P) \in M(P) + \sum_{Q: Q<P} t^{-1}\Z[t^{-1}] M(Q).
\end{equation*}

\textup{(2)} The specialization of $L(P)$ at $t=1$ coincides with the
simple module with Drinfeld polynomial $P$.
\end{Theorem}

As we mentioned in the introduction, the relation between $M(P)$ and
$L(P)$ in $\bfR_t$ (not in its specialization) can be understood by
a Jantzen filtration \cite{Gr}.

For a later purpose we define matrices in the Laurent polynomial ring
of $t$:
\begin{gather*}
  c_{PQ}(t) \defeq
  \text{the coefficient of $e^Q$ in $\chi_{\ve,t}(M(P))$},
\\
  (c^{PQ}(t)) \defeq (c_{PQ}(t))^{-1},
\\
  M(P) = \sum_Q Z_{PQ}(t) L(Q).
\end{gather*}

When $\ve$ is not a root of unity, there is an isomorphism between
$\bfR_t$ and the dual of the Grothendieck group of a category of
perverse sheaves on affine graded quiver varieties in
\cite[\S14]{Na-qaff}. And the full detailed proof of the above theorem
was already explained in \cite{Na-ann}.
However, the latter group becomes larger when $\ve$ is a root of
unity.
So we modify $\bfR_t$ to $\widetilde\bfR_t$, and give a proof of the
above theorem in this $\widetilde\bfR_t$.

Let us fix an $I$-tuple of polynomials $P$ throughout this section.
Let $\mathcal I$ be the set of {\it l\/}--dominant monomials
$m\in\widehat{\mathscr Y}_t$ such that $m \le e^P$. We consider
$\Z[t,t^{-1}]$-module with basis $\mathcal I$, and denote it by
$\widetilde\bfR_t$.

For each monomial $m\in\mathcal I$, let $P_m$ be an $I$-tuple of
polynomials given by
\(
  (P_m)_i(u) \defeq \prod_a (1 - ua)^{u_{i,a}(m)}
\).
In other words, $P_m$ is determined so that $\widehat\Pi(e^{P_m}) =
\widehat\Pi(m)$.
If $\ve$ is not a root of unity, then $P_m = P_{m'}$ implies $m = m'$
by the invertibility of the $\ve$-analog of the Cartan matrix. But it
is not true in general. This is the reason why we need a modification.

We modify $\widehat\chi_{\ve,t}$ of the standard module $M(P_m)$ so
that it has the image in $\widetilde\bfR_t$ as follows: If
\begin{equation*}
   \widehat\chi_{\ve,t}(M(P_m)) = \sum_{n} a_{n,m}(t) e^{P_m} n,
\end{equation*}
then, we define
\begin{equation*}
   M_m \defeq 
    \sum_{n^*} 
     a_{n^*,m}(t)\, t^{-d(e^{P_m}n^*, e^{P_m}n^*)} m n^*,
\end{equation*}
where the summation runs {\it only\/} over $n^*$ such that $mn^*$ is {\it
  l\/}--dominant.
The $M_m$ is contained in $\widetilde\bfR_t$ by Axiom~1. And
$\widehat\chi_{\ve,t}(M(P_m))$ is recoved from $M_m$.
Let us denote the coefficient of $mn^{*}$ by $c_{mn}(t)$, where
$n = mn^{*}\in\mathcal I$, that is
\begin{equation}\label{eq:c}
   M_m = \sum_n c_{mn}(t) n.
\end{equation}
We have $c_{mm} = 1$ and $c_{mn}(t) = 0$ for $n\nleq m$. In
particular, $\{ M_m \}_m$ is a base of $\widetilde\bfR_t$.

We define an involution
$\setbox5=\hbox{A}\overline{\rule{0mm}{\ht5}\hspace*{\wd5}}$ on
$\widetilde\bfR_t$ by
\begin{equation*}
  \overline{t} = t^{-1}, \quad \overline{m} = m.
\end{equation*}

We define a map $\widetilde\bfR_t\to\bfR_t$ by $M_m\mapsto M(P_m)$.
When $\ve$ is not a root of unity, this map is injective and the image
is the submodule spanned by $M(P')$'s such that
$\widehat\Pi(e^{P'})\le\widehat\Pi(e^P)$.
The map intertwines the involutions.

We have
\begin{equation*}
   \overline{M_m} = \sum_n c_{mn}(t^{-1}) n
   = \sum_{n,s} c_{mn}(t^{-1})c^{ns}(t) M_s,
\end{equation*}
where $(c^{ns}(t))$ is the inverse matrix of $(c_{mn}(t))$.
Let
\begin{equation}\label{eq:u}
   u_{mn}(t) \defeq \sum_s c_{ms}(t^{-1}) c^{sn}(t),
   \qquad\text{or equivalently }
   \overline{M_m} = \sum_n u_{mn}(t) M_n.
\end{equation}
By axioms, $u_{mm}(t) = 1$ and $u_{mn}(t) = 0$ if $n\nleq m$.

\begin{Lemma}
There exists a unique element $L_m\in\widetilde\bfR_t$ such that
\begin{equation*}
\overline{L_m} = L_m, \qquad
    L_m \in M_m + \sum_{n: n<m} t^{-1}\Z[t^{-1}] M_n.
\end{equation*}
\end{Lemma}

Although the proof is exactly the same as \cite[7.10]{Lu:can}, we give 
it for the sake of the reader.
\begin{proof}
Let
\[
    M_m = \sum_{n\le m} Z_{mn}(t) L_n.
\]
Then the condition for $\{ L_m\}$ is equivalent to the following
system:
\begin{subequations}\label{eq:ZZ}
\begin{align}
   & Z_{mm}(t) = 1, \quad
   \text{$Z_{mn}(t)\in t^{-1}\Z[t^{-1}]$ for $n < m$}, \label{eq:Z1}
\\
   & Z_{mn}(t^{-1}) = \sum_{s: n \le s\le m} u_{ms}(t) Z_{sn}(t).
   \label{eq:Z2}
\end{align}
\end{subequations}
The equation can be rewritten as
\begin{equation*}
   Z_{mn}(t^{-1}) - Z_{mn}(t)
   = \sum_{s: n \le s < m} u_{ms}(t) Z_{sn}(t).
\end{equation*}
Let $F_{mn}(t)$ be the right hand side.
We can solve this system uniquely by induction: If $Z_{sn}(t)$'s are
given, $Z_{mn}(t)$ is uniquely determined by the above equation and
$Z_{mn}(t)\in t^{-1}\Z[t^{-1}]$, provided $F_{mn}(t^{-1}) = -
F_{mn}(t)$. We can check this condition by the induction hypothesis:
\begin{equation*}
\begin{split}
   F_{mn}(t^{-1}) & =
    \sum_{s: n\le s < m} u_{ms}(t^{-1}) Z_{sn}(t^{-1})
\\
   & = \sum_{s: n\le s < m} \sum_{t:n\le t\le s}
     u_{ms}(t^{-1}) u_{st}(t) Z_{tn}(t)
\\
   & = - \sum_{t: n\le t < m} u_{mt}(t) Z_{tn}(t)
     = - F_{mn}(t),
\end{split}
\end{equation*}
where we have used
\(
  \sum_{s:t\le s\le m} u_{ms}(t^{-1}) u_{st}(t) = 0
\)
for $t < m$.
\end{proof}

The proof of \thmref{thm:character}(1) is exactly the same. Since the
map $\widetilde\bfR_t\to\bfR_t$ intertwines the involution, the image
of $L_m$ is equal to $L(P_m)$. Therefore \thmref{thm:character}(2) is
equivalent to the following statement:

\begin{Theorem}\label{thm:character2}
The multiplicity $[M(P):L(Q)]$ is equal to
\begin{equation*}
    \sum_n Z_{e^P, n}(1),
\end{equation*}
where the summation is over the set $\{ n\mid e^Q = \widehat\Pi(n) \}$.
\end{Theorem}

The following proof is just a modification of that given in
\cite{Na-ann}.

We choose $W$ so that $e^P = e^W$ as before. 
Let $D^b(\N_0(\infty,W))$ be the bounded derived category of complexes
of sheaves whose cohomology sheaves are constant along each connected
component of a stratum $\Nreg(V,W)$ of \eqref{eq:stratum}. 
(The connectedness of $\Nreg(V,W)$ is {\it not\/} known.)
If $\Nreg(V,W)^\alpha$ is a connected component of $\Nreg(V,W)$, then
$IC(\Nreg(V,W)^\alpha)$ be the intersection homology complex
associated with the constant local system $\C_{\Nreg(V,W)^\alpha}$ on
$\Nreg(V,W)^\alpha$.
Then $D^b(\N_0(\infty,W))$ is the category of complex of sheaves which
are finite direct sums of complexes of the forms
$IC(\Nreg(V,W)^\alpha)[d]$ for various $V$, $\alpha$ and $d\in\Z$,
thanks to the existence of transversal slices \cite[\S3]{Na-qaff}.


We associate a monomial $m = e^V e^W$ to each $[V]$. It gives us a
bijective correspondence between the set of monomials $m$ with $m\le
e^P$ and the set of isomorphism classes of $I\times\C^*$-graded vector
spaces. If $\Nreg(V,W)\neq\emptyset$, the corresponding monomial $m$
is {\it l\/}--dominant, i.e., $m\in\mathcal I$.
We choose a point in $\Nreg(V,W)$ and denote it by $x_m$.

Let $\C_{\N(V',W)}$ be the constant local system on $\N(V',W)$.
Then $\pi_*\C_{\N(V',W)}$ is an object of $D^b(\N_0(\infty,W))$ again
by the transversal slice argument. From the decomposition theorem of
Beilinson-Bernstein-Deligne, we have
\begin{equation}\label{eq:decomp}
   \pi_*(\C_{\N(V',W)}[\dim_\C \N(V',W)]) 
        \cong \bigoplus_{V,\alpha,k}
      L_{V,\alpha,k}(V',W)\otimes IC(\Nreg(V,W)^\alpha)[k]
\end{equation}
for some vector space $L_{V',\alpha,k}(V,W)$ \cite[14.3.2]{Na-qaff}.
We set
\begin{equation}\label{eq:L}
  L_{mn}(t) \defeq \sum_k \dim L_{V,\alpha,k}(V',W)\, t^{-k},
\end{equation}
where $V$, $V'$ are determined so that $m = e^{V}e^{W}$, $n =
e^{V'}e^W$. By the description of the transversal slice
\cite[\S3]{Na-qaff}, $\dim L_{V,\alpha,k}$ is independent of $\alpha$.
So $\alpha$ can disappear in the left hand side.
Applying the Verdier duality to the both hand side of
\eqref{eq:decomp} and using the self-duality of
$\pi_*(\C_{\N(V',W)}[\dim_\C\N(V',W)])$ and $IC(\Nreg(V,W)^\alpha)$,
we find $L_{m'm}(t) = L_{m'm}(t^{-1})$.

By our definition of $M_m$, we have
\begin{equation*}
   M_m = \sum_{[V_n]} t^{-\dim\N(V_n,W_m)} P_t(\NLa(V_n,W_m))\, m e^{V_n},
\end{equation*}
where $W_m$ is given by $\dim (W_m)_{i,a} = u_{i,a}(m)$.
By \cite[\S3]{Na-qaff}, this is equal to
\begin{equation}\label{eq:C_RQ}
\begin{split}
   M_m & = \sum_{[V_n]} t^{-\dim\N(V_n,W_m)}
   P_t\left(\pi^{-1}(x_m)\cap\N(V_n\oplus V_m,W)\right)\, m e^{V_n}
\\
   & = \sum_{[V']} \sum_k t^{\dim\N(V_m,W)-k}
   \dim H^k(i_{x_m}^! \pi_* \C_{\N(V',W)}[\dim \N(V',W)])
     \, e^{V'}e^{W},
\end{split}
\end{equation}
where $V_m$ is given so that $e^{V_m} e^{W} = m$.

Let
\begin{equation*}
   Z_{mn}(t) \defeq
   \sum_{k,\alpha}
    \dim H^k(i_{x_m}^! IC(\Nreg(V,W)^\alpha))\, t^{\dim\Nreg(V_m,W)-k},
\end{equation*}
where $n = e^V e^W$. By the defining property of the intersection
homology, we have \eqref{eq:Z1} and
$Z_{mn}(t) = 0$ if $n \nleq m$.

Substituting \eqref{eq:decomp} into \eqref{eq:C_RQ}, we get
\begin{equation}\label{eq:C=ZL}
   c_{mn}(t) = \sum_{s} Z_{ms}(t) L_{sn}(t).
\end{equation}
Now $L_{sn}(t) = L_{sn}(t^{-1})$ and \eqref{eq:u} imply \eqref{eq:Z2}.

Let $Z^\bullet(W)$ be the fiber product
$\N(W)\times_{\N_0(\infty,W)}\N(W)$.
Let $\mathcal A = H_*(Z^\bullet(W),\C)$ be its Borel-Moore homology
group, equipped with an algebra structure by the convolution (see
\cite[14.2]{Na-qaff}).
Taking direct sum with respect to $V'$ in \eqref{eq:decomp}, we have a
linear isomorphism (forgetting gradings)
\begin{equation*}
    \pi_*(\C_{\N(W)}) =
    \bigoplus_{V,\alpha} L_{V,\alpha}\otimes IC(\Nreg(V,W)^\alpha),
\end{equation*}
where $L_{V,\alpha} = \bigoplus_{[V'],k} L_{V,\alpha,k}(V',W)$.
By a general theory (see \cite{Gi-book} or \cite[14.2]{Na-qaff}), $\{
L_{V,\alpha}\}$ is a complete set of mutually nonisomorphic simple
$\mathcal A$-modules.
Moreover, taking $H^*(i_{x_m}^!\ )$ of both hand sides, we have
\begin{equation*}
    H(\pi^{-1}(x_m),\C)) =
    \bigoplus_{V,\alpha}
     L_{V,\alpha}\otimes H^*(i_{x_m}^!IC(\Nreg(V,W)^\alpha)),
\end{equation*}
which is an equality in the Grothendieck group of $\mathcal A$-modules.
Here the $\mathcal A$-module structure on the right hand side is given by
$a:\xi\otimes\xi'\mapsto a\xi\otimes\xi'$.

By \cite[\S13]{Na-qaff}, there exists an algebra homomorphism
$\Ule\to\mathcal A$. Moreover \cite[\S14.3]{Na-qaff}, each
$L_{V,\alpha}$ is a simple {\it l\/}--highest weight $\Ule$-module.
Its Drinfeld polynomial is $Q$ such that $\widehat\Pi(e^V e^W) = e^Q$.
(It is possible to have two different $V$, $V'$ give isomorphic
$\Ule$-modules.)
Combining with above discussions, we get \thmref{thm:character2}.

\begin{Remark}
If one enlarges the commutative subalgebra $\Ule^0$ of $\Ule$, then
he/she can recover a bijective correspondence between simple
$\Ule$-modules, and strata of affine quiver varieties. When $\g$ is of
type $A_n$, such the enlargement is ${\mathbf
U}_{\varepsilon}({\mathbf L}\mathfrak{gl}_{n+1})$.  (cf. \cite{FM3})
\end{Remark}

\section{The specialization at $\ve = \pm 1$}\label{sec:pm 1}

When $\ve = \pm 1$, simple modules can be described explicitly
\cite[\S4.8]{FM2}. We study their $\widehat\chi_{\ve,t}$ in this
section.

Let $P$ be an $I$-tuple of polynomials. We choose $I\times\C^*$-graded
vector space $W$ so that $e^W = e^P$ as before.

First consider the case $\ve = 1$. The $W$ can be considered as a
collection of $I$-graded vector spaces $\{ W^a\}_{a\in\C^*}$, where
$W^a_i \defeq W_{i}(a)$. Then from the definition of cyclic quiver
varieties, it is clear that we have
\begin{equation*}
   \N(W) \cong \prod_a \M(W^a),\qquad
   \N_0(\infty,W) \cong \prod_a \M_0(\infty,W^a),
\end{equation*}
Here $\M(W^a)$, $\M_0(\infty,W^a)$ are the original quiver variety
corresponding to $W^a$.
Let $P^a$ be an $I$-tuple of polynomial defined by
$P^a_i(u) = (1 - au)^{\dim W_i(a)}$. The $P^a$ is, of course,
determined directly from $P$. From above description, we have
\begin{equation}\label{eq:prod}
   M(P) = \bigotimes_a M(P^a), \qquad
   \widehat\chi_{\ve,t}(M(P)) = \prod_a \widehat\chi_{\ve,t}(M(P^a)).
\end{equation}
The latter also follows directly from Axiom~3.

Next consider the case $\ve = -1$. We choose and fix a funcion
$o\colon I\to \{\pm 1\}$ such that $o(i) = -o(j)$ if $a_{ij}\neq 0$,
$i\neq j$. We define an $I$-graded vector space $W^a$ by $W^a_i \defeq
W_i(o(i)a)$. Then we have
\begin{equation*}
   \N(W) \cong \prod_a \M(W^a),\qquad
   \N_0(\infty,W) \cong \prod_a \M_0(\infty,W^a).
\end{equation*}
More precisely, $\N(V,W) = \prod_a \M(V^a,W^a)$ with $V^a \defeq
\bigoplus_i V_i(-o(i)a)$.
Let $P^a$ be an $I$-tuple of polynomial defined by
$P^a_i(u) = (1 - o(i)au)^{\dim W_i(o(i)a)}$. The $P^a$ is again
determined directly from $P$. We have \eqref{eq:prod} also in this case.

Recall that we have an algebra homomorphism ${\mathbf U}_{\ve}(\g)\to
\Ule$ \eqref{eq:subalg}. By \cite[\S33]{Lu-book}, ${\mathbf
U}_{-1}(\g)$ is isomorphic to ${\mathbf U}_{1}(\g)$. Moreover, the
universal enveloping algebra ${\mathbf U}(\g)$ of $\g$ is isomorphic
to the quotient of ${\mathbf U}_1(\g)$ by the ideal generated by $q^h -
1$ ($h\in P^*$) \cite[9.3.10]{CP-book}. In particular, the category of
type $1$ finite dimensional ${\mathbf U}_{\ve}(\g)$-modules is
equivalent to the category of finite dimensional $\g$-modules.
Therefore we consider $\operatorname{Res}M(P)$ as a $\g$-module.

Thanks to the fact that $\pi\colon\M(W)\to\M_0(\infty,W)$ is
semismall, we have the following \cite[\S15]{Na-qaff}:
\begin{Theorem}
\textup{(1)} 
\(
    L(P) = \bigotimes_a L(P^a).
\)

\textup{(2)} For each $a$, $\operatorname{Res}(L(P^a))$ is simple as a
$\g$-module. Its highest wight is $\Lambda^a = \sum_i \deg P^a_i
\Lambda_i$
\end{Theorem}

We want to interpret this result from $\widehat\chi_{\ve,t}$. We
may assume that there exists only one nontrivial $P^a$. All other
$P^b$'s are $1$.

We identify an $I$-tuple of polynomials $P$ whose roots are $1/a$ with
a dominant weight by
\begin{equation*}
    P \mapsto \deg P \defeq \sum_i \deg P_i\, \Lambda_i, \qquad
    \lambda = \sum_i \lambda_i\Lambda_i \mapsto 
     P_\lambda; (P_\lambda)_i(u) = (1 - au)^{\lambda_i}.
\end{equation*}

We give an explicit formula of $Z_{PQ}(t)$, not based on
inductive procedure:
\begin{Theorem}
  \textup{(1)} $Z_{PQ}(1)$ is equal to the multiplicity of the
  simple $\g$-module $L(\deg Q)$ of highest weight $\deg Q$ in
  $\operatorname{Res}M(P)$.
  
  \textup{(2)} $c_{PQ}(t)$ is a polynomial in $t^{-1}$, so
  $c_{PQ}(\infty)$ makes sense.
  
  \textup{(3)} We have
\begin{equation*}
  \begin{split}
   & \chi_{\ve,t}(L(P)) = \sum_Q c_{PQ}(\infty)\, e^Q
   + \text{non {\it l\/}--dominant terms},
\\
   &\qquad\text{or equivalently }
   Z_{PQ}(t) = \sum_R c_{PR}(t)c^{RQ}(\infty).
  \end{split}
\end{equation*}

\textup{(4)} The coefficient $c_{PQ}(\infty)$ is equal to the
weight multiplicity of the dominant weight $\deg Q$ in $L(\deg P)$.
\end{Theorem}

\begin{proof}
(1) is clear.

(2),(3) By the fact that $\pi\colon\M(W^a)\to\M_0(\infty,W^a)$ is
semismall, we have $L_{V,\alpha,k} = 0$ (in \eqref{eq:decomp}) for
$k\neq 0$, hence $L_{PQ}(t)$ (in \eqref{eq:L}) is a constant. Then
$c_{PQ}(t)$ is a polynomial in $t^{-1}$ by \eqref{eq:C=ZL} and
\eqref{eq:Z1}.  Therefore $c_{PQ}(\infty)$ makes sense. We have
$c_{PQ}(\infty) = L_{PQ}(t) = L_{PQ}(0)$ again by \eqref{eq:C=ZL} and
\eqref{eq:Z1}. Thus we get the assertion.

(4) By (3), we have
\begin{equation*}
   \chi(L(P)) = \sum_Q c_{PQ}(\infty)e^{\deg Q} + \text{non dominant terms}.
\end{equation*}
Since $\chi$ is the ordinary character, the assertion is clear.
\end{proof}

Note that the multiplicity of a simple (resp.\ Weyl) ${\mathbf
U}_\ve(\g)$-module in $\operatorname{Res}M(P)$ for generic (resp.\ a
root of unity) $\ve$ is independent of $\ve$. Therefore $Z_{PQ}(1)$
gives it. (See \cite[\S7]{Na-ann}.)

\section{Conjecture}\label{sec:conjecture}

There is a large amount of literature on finite dimensional
$\Ule$-modules.
Some special classes of simple finite dimensional $\Ule$-modules are
studied intensively: tame modules \cite{NT} and Kirillov-Reshetikhin
modules \cite{HKOTY} (see also the references therein). For tame
modules, there are explicit formulae of $\chi_\ve$ in terms of Young
tableaux. For Kirillov-Reshetikhin modules, there are conjectural
explicit formula of $\chi$ (i.e., decomposition numbers of
restrictions to ${\mathbf U}_\ve(\g)$-modules).

Although our computation applies to {\it arbitrary\/} simple modules,
our polynomials $Z_{PQ}(t)$ are determined recursively, and it is
difficult to obtain explicit formulae in general. Thus those modules
should have a very special feature among arbitrary modules.
For Kazhdan-Lusztig polynomials, a special classes is known to have
explicit formulae. Those are Kazhdan-Lusztig polynomials for
Grassmannians studied \cite{LS}. A geometric interpretation was given
in \cite{Ze}. Based on an analogy between Kazhdan-Lusztig polynomials
and our polynomials, we propose a class of finite dimensional
$\Ule$-modules. It is a class of {\it small\/} standard modules.

\begin{Definition}
(1) A finite dimensional $\Ule$-module $M$ is said {\it special\/} if
it satisfies the condition in \thmref{thm:cons}(2),
i.e., $\chi_\ve(M)$ contains only one {\it l\/}--dominant monomial.

(2) Let $M(P)$ be a standard module with {\it l\/}--highest weight
$P$. We say $M(P)$ is {\it small\/} if $c_{QR}(t)\in t^{-1}\Z[t^{-1}]$
for any $Q$, $R \le P$ with $Q\neq R$.
Similarly $M(P)$ is called {\it semismall\/} if
$c_{QR}(t)\in \Z[t^{-1}]$ for any $Q,R\le P$.
\end{Definition}

\begin{Remark}
(1) By the geometric definition of $\widehat\chi_{\ve,t}$
\eqref{eq:geomdef}, $M(P)$ is (semi)small if and only if
\(
   \pi\colon\N(V,W) \to \N_0(V,W)
\)
is (semi)small for any $V$ such that $e^V e^W$ is {\it l\/}--dominant.

(2) By definition, $M(Q)$ is (semi)small if $M(P)$ is (semi)small and
$Q\le P$.

(3) The (semi)smallness of $M(P)$ is related to (semi)tightness of
monomials in $\mathbf U^-_q$ \cite{Lu:tight}.
\end{Remark}

Since a finite dimensional simple $\Ule$-module contains at least one
{\it l\/}--dominant monomial, namely the one corresponding to the {\it 
l\/}--highest weight vector, a special module is automatically
simple. The converse is {\it not\/} true in general. For example, if
$\g = \algsl_2$, $P = (1-u)^2(1-\ve^2 u)$, then one can compute (say,
by our algorithm)
\[
  \chi_\ve(L(P)) = Y_{1,1}^2 Y_{1,\ve^2} + Y_{1,1} 
  + Y_{1,1}^2 Y_{1,\ve^4}^{-1} 
  + 2 Y_{1,1}Y_{1,\ve^2}^{-1}Y_{1,\ve^4}^{-1}
  + Y_{1,\ve^2}^{-2}Y_{1,\ve^4}^{-1}.
\]
This has two {\it l\/}-dominant monomial terms.

\begin{Theorem}
  Suppose $M(P)$ is small. Then for any $I$-tuple of polynomials $Q\le
  P$, corresponding simple module $L(Q)$ is special.
\end{Theorem}

\begin{proof}
By the characterization of $Z_{QR}(t)$ in \eqref{eq:ZZ}, we have
$Z_{QR}(t) = c_{QR}(t)$ for all $Q,R\le P$. Therefore
\begin{equation*}
   \sum_R Z_{QR}(t) \chi_{\ve,t}(L(R)) = \chi_{\ve,t}(M(Q))
   = \sum_R Z_{QR}(t) e^R + \text{non {\it l\/}--dominant terms}.
\end{equation*}
Hence we have $\chi_{\ve,t}(L(R))$ is $e^R$ plus non {\it
  l\/}--dominant terms.
\end{proof}

\begin{Conjecture}
  Standard modules corresponding to tame modules and
  Kirillov-Reshetikhin modules are small.
\end{Conjecture}

\end{document}